\documentclass[11pt, reqno]{amsart}

 \usepackage{amsmath}
 \usepackage{amssymb}
\usepackage{bm}

\DeclareMathAccent{\mathring}{\mathalpha}{operators}{"17}

 \usepackage{color}

\newcommand{\mysection}[1]{\section{#1}
      \setcounter{equation}{0}}

\newtheorem{theorem}{Theorem}[section]
\newtheorem{lemma}[theorem]{Lemma}

\newtheorem{corollary}[theorem]{Corollary} 

\theoremstyle{definition}
\newtheorem{assumption}{Assumption}[section]
\newtheorem{definition}{Definition}[section]
\theoremstyle{remark}
\newtheorem{remark}{Remark}[section]

\newcommand{\tr}{\text{\rm tr}\,}
\newcommand{\loc}{\text{\rm loc}}

 \makeatletter 
 \def\dashint{%
 \operatorname%
 {\,\,\text{\bf--}\kern-.98em\DOTSI\intop\ilimits@\!\!}}
\def\ninf{\qopname\relax\@empty{inf\phantom{p}\!\!\!}}
 \makeatother

\newcommand\bbeta{\text{\raise-.2ex\hbox{$\bm{\beta}$}}}

\newcommand\esssup{\operatornamewithlimits{esssup}}

\newcommand\bR{\mathbb{R}}

\newcommand\bS{\mathbb{S}}

\newcommand\cB{\mathcal{B}}

\newcommand\cF{\mathcal{F}}

\newcommand\cM{\mathcal{M}}

\newcommand\cN{\mathcal{N}}

\newcommand\dist{{\rm dist}\,}

\newcommand\osc{\operatornamewithlimits{osc}}

\begin{document}

\title[Diffusions with drift in $L_{d}$]
{On diffusion processes with drift in $L_{d}$}

\author{N.V. Krylov}
 
\email{nkrylov@umn.edu}
\address{127 Vincent Hall, University of Minnesota,
 Minneapolis, MN, 55455}



\begin{abstract}
We investigate properties of Markov quasi-diffusion processes
corresponding to elliptic operators $L=a^{ij}D_{ij}+b^{i}D_{i}$,
acting on functions on $\mathbb{R}^{d}$, with 
measurable coefficients, bounded 
and uniformly elliptic $a$ and $b\in L_{d}(\mathbb{R}^{d})$. We 
show that each of them is strong Markov with 
strong Feller transition
semigroup $T_{t}$, which is also a continuous 
bounded semigroup in
$L_{d_{0}}(\mathbb{R}^{d})$ for some $d_{0}\in (d/2, d)$. 
We show 
that  $T_{t}$, $t>0$, has a kernel $p_{t}(x,y)$ which is 
summable in $y$ to the power of
$d_{0}/(d_{0}-1)$. This leads to the parabolic 
Aleksandrov estimate
with power of summability $d_{0}$ instead of the 
usual $d+1$. For the probabilistic
solution, associated with such a process, 
of the problem $Lu=f$
in a bounded domain $D\subset\mathbb{R}^{d}$ 
with boundary condition
$u=g$, where $f\in L_{d_{0}}(D)$ and $g$ is 
bounded, we show that it is
H\"older continuous. Parabolic version of this
 problem is treated as well. 
We also prove Harnack's inequality for harmonic 
and caloric functions
associated with such a process. Finally, we show that 
the probabilistic solutions
are $L_{d_{0}}$-viscosity solutions.
  
\end{abstract}

\maketitle

\mysection{Introduction}
                                  \label{section 12.29.1}

Let $\bR^{d}$ be a Euclidean space of
points $x=(x^{1},...,x^{d})$. For a fixed 
throughout the article $\delta\in(0,1)$ define
$\bS_{\delta}$ as the set of $d\times d$ symmetric
matrices whose eigenvalues are between $\delta$
and $\delta^{-1}$. Fix a constant $\|b\|\in(0,\infty)$.
In this article we consider and discuss
only uniformly nondegenerate
processes with bounded diffusion coefficient.

\begin{assumption}
                                \label{assumption 12.29.1}
We are given a Borel measurable
$\bS_{\delta}$-valued function $a=a(x)$ and
a Borel measurable
$\bR^{d}$-valued function $b=b(x)$ such that
$$
\|b\|_{L_{d}(\bR^{d})}\leq \|b\|.
$$
\end{assumption}
Define
\begin{equation}
                                         \label{1.12.2}
D_{i}=\frac{\partial}{\partial x^{i}},
\quad D_{ij}=D_{i}D_{j},\quad 
L=(1/2)a^{ij}(x)D_{ij}+b^{i}(x)D_{i}.
\end{equation}

The goal of this article is to investigate 
(time-homogeneous Markov)
quasi-diffusion processes corresponding
to $L$. In the more modern terminology from
\cite{SV_79} these are called diffusion
processes, but at this point and later on
we will follow the terminology from
\cite{Dy_63} in which
the notion of diffusion processes is defined
differently from
\cite{SV_79}.

The definition of time-homogeneous diffusion processes
 first appeared
in the book by Dynkin in 1963, \cite{Dy_63}, where he also constructs
 diffusion  processes corresponding to elliptic
 operators as in \eqref{1.12.2}
with bounded and H\"older continuous coefficients,
such that the matrix $ (a_{ij}(x))$ is uniformly strictly
positive.

If $ x_{t}(x)$, $ t\ge 0$, is a family of continuous processes
on $ \bR^{d}$, parametrized by $ x\in \bR^{d}$, and
 the family is a diffusion process corresponding to the above $ L$
in Dynkin's sense, then, for any bounded domain $ D\subset\bR^{d}$
and smooth function $ u$,
$$ 
u(x)=E_{x}\bigg[u(x_{\tau_{D}})-\int_{0}^{\tau_{D}}Lu(x_{t})\,dt
 \bigg],
$$
where $ \tau_{D}=\tau_{D}(x)$ is the first exit time
of $ x_{t}(x)$ from $ D$.

 The author took the above property as the definition of
{\em quasi-diffusion\/}
 process and in 1966 constructed such process 
under the assumptions that the matrix 
$ (a_{ij}(x))$ is uniformly strictly
positive, is continuous and $ b$ is Borel bounded.
The domain of definition
of the corresponding generator of the constructed process
was also described,
which provides the so-called weak uniqueness of the process
with this generator. Later on in 1973, (\cite{Kr_73}), when It\^o's
formula was extended to $ W^{2}_{p}$ functions,
it became obvious that the quasi-diffusion
processes corresponding to the operators satisfying
the above mentioned condition
are weakly unique.
 
Two years earlier Tanaka constructed Dynkin's diffusion
processes when $ b$ is also continuous.
No uniqueness was implied in his paper.

Quasi-diffusion processes are characterized by the property that,
for any smooth function $ u(x)$ and starting point $ x$, the process
$$ 
u(x_{t}(x))-\int_{0}^{t}Lu(x_{s}(x))\,ds
$$
is a local martingale. Stroock and Varadhan (1969),
\cite{SV_79}, took the time-inhomogeneous version
of this property as the definition
of diffusion process and proved existence and weak
uniqueness under the condition that $  a,b$ are bounded,
$ a$ is uniformly continuous in $  x$ uniformly in $  t$,
and $  a$ is uniformly nondegenerate. The proof
of uniqueness is based on 
the solvability of parabolic equations with
coefficients depending only on time, a result borrowed from 
PDE, and
the estimate
\begin{equation} 
                                               \label{1.2.2}
E_{s,x}\int_{0}^{1}f(s+t,x_{t})\,dt\leq
N\bigg(\int_{0}^{1}\int_{\bR^{d}}|f(s+t,y)|^{p}\,dtdx\bigg)^{1/p},
\end{equation}
which is achieved by a quite clever argument.
This argument, however, is heavily based on
the uniform continuity of $ a$ with respect to $ x$.

In 1974 \cite{Kr_74}
the author proved that estimate \eqref{1.2.2} holds true
not only for solutions of stochastic equations but also 
in the case that $  a,b$  
are any progressively measurable
bounded functions such that $  a$ is uniformly nondegenerate.
The method of proof is different from the one used by 
Stroock and Varadhan (and the range of $p$ is more restrictive).

 Recall that the first quasi-diffusion strong Markov
 processes with
{\em bounded\/} Borel $b$ and Borel uniformly nondegenerate $a$ were
constructed in \cite{Kr_73}. This construction was carried
over to the case of time-inhomogeneous processes with jumps in
\cite{AP_77}. A different approach  again 
when $b$ is bounded,
 based
on Krylov-Safonov estimates, is carried out
in \cite{Ba_98} and produced a 
particular strong Markov process
with strong Feller resolvent.

 It is worth mentioning that
with sufficiently regular diffusion  
matrix $a$  and the drift which is a
generalized function of the type of the derivative of a
measure generalized diffusion processes are constructed
in \cite{Po_82}.  The case of
 time-dependent regular $a$  and $b$
summable to powers which are  different in $t$ and $x$
attracted a very extensive attention. In that regard
the reader can consult
\cite{ZZ_20}, \cite{Po_82}, and references therein.
Our main emphasis here is on 
time-homogeneous
$b\in L_{d}(\bR^{d})$
and Borel $\bS_{\delta}$-valued $a$.

In case $b\in L_{d}$
we established in \cite{Kr_19_1} the existence of a strong
Markov quasi-diffusion process corresponding to $L$.
Our goal in this article is to investigate
properties of just
Markov quasi-diffusion processes corresponding to $L$
 regardless of how
they appeared or were constructed.

In particular, we show that each of them is
strong Markov with strong Feller transition
semigroup $T_{t}$, which is also a continuous bounded semigroup
in $L_{p}(\bR^{d})$ for $p\in[d_{0},d)$,
where $d_{0}\in(d/2,d)$.
Our estimate \eqref{12.2.1} implies that $T_{t}$
has a kernel $p_{t}(x,y)$ which is summable
with respect to  $y$ to the power of $d_{0}/(d_{0}-1)$.
This leads to the parabolic Aleksandrov estimate
with power of summability $d_{0}$ instead of
usual $d+1$ (see Corollary \ref{corollary 12.5.3} 
for the probabilistic version and
Theorem \ref{theorem 12.30.1} for the PDE version).

For the probabilistic solution,
associated with such a process, of the problem
$Lu+f=0$ in a bounded domain $D\subset \bR^{d}$
with boundary condition $u=g$, where $f\in L_{d_{0}}(D)$
and $g$ is Borel bounded, we show that it is H\"older
continuous in $D$. Parabolic version of this
problem is treated as well. We also prove Harnack's inequality for
harmonic and caloric functions associated 
with such a process. Finally we show that probabilistic
solutions of the corresponding elliptic
equations are $L_{d_{0}}$-viscosity solutions.

In our arguments self-similar transformations
play a very important role. Observe that under such
transformations $a$ and $b$ change, but the new $a$
is still in $\bS_{\delta}$ and the $L_{d}(\bR^{d})$-norm
of the new $b$ is majorated by {\em the same\/} 
number $\|b\|$.

We  finish the introduction with some notation.
For $T,R>0$, $(t,x)\in \bR^{d+1}=\{(t,x):t\in\bR,x\in\bR^{d}\}$ define 
$$
B_{R}(x)=\{y\in\bR^{d}:|y-x|<R\},\quad
B_{R} =B_{R}(0), \quad C _{T,R}=[0,T)\times B_{R},
$$ 
$$
\quad C_{T,R}(t,x)=(t,x)+C _{T,R},
\quad C _{R}(t,x)=C _{R^{2},R}(t,x),\quad C _{R} =C _{R}(0,0),
  $$
$$
\partial_{t}=\partial/\partial t.
$$

In the proofs of various results  we use
the symbol $N$  to denote finite 
nonnegative constants
which may change from one occurrence to another and
we do not always specify on which data these  constants
depend. In these cases the reader should remember
that, if in the statement of a result there are constants
called $N$ which are claimed to depend only on certain
parameters, then in the proof of the result
the constants $N$ also depend only on the same
parameters unless specifically stated otherwise.
Of course, if we write 
$
N=N(...),
$
 this means that $N$ depends only
on what is inside the parentheses. 
Another point is that when we say that certain constants
depend only on such and such parameters we mean, in particular,
  that the dependence is such that these constants stay bounded
as the parameters vary in compact subsets of their ranges.

\mysection{Diffusions and It\^o stochastic equations}
                                      \label{section 12.31.1}

Suppose that we are given a 
a quasi-diffusion process corresponding to $L$, that is,
we are given a
continuous Markov process $X=(x_{t},\infty,\cM_{t},P_{x})$
(the terminology taken from \cite{Dy_63})
in $\bR^{d}$ such that 
  for any $x\in\bR^{d}$ and $t\geq0$
\begin{equation}
                                               \label{11.4.3}
E_{x} \int_{0}^{t}|b(x_{s})|\,ds<\infty 
\end{equation}
and for any
twice continuously differentiable function $u$
with compact support
$$
u(x)=E_{x}u(x_{t})-E_{x}\int_{0}^{t}Lu(x_{s})\,ds.
$$

\begin{remark}
                                \label{remark 11.4.1}
These requirements would be unrealistic if $b$
were only in $L_{p}$ with $p<d$ (see \cite{Kr_19_1}).
\end{remark}
 
Define
$\tau_{R} $ as the first exit time of $ x_{t}$ from  
$B_{R}$ (equal to infinity if $ x_{t}$ never exits from $B_{ R}$).
  {\em This notation $\tau_{R}$ is used throughout
the article}.

Denote by $\cN_{t}$ the $\sigma$-field in $\Omega$
generated by the events $\{\omega:x_{s}\in\Gamma\}$
for all $s\leq t$ and Borel $\Gamma\subset\bR^{d}$
and let $\bar \cN^{x_{0}}_{t}$ be the completion
of $\cN_{t}$ with respect to $(\cN_{\infty},P_{x_{0}})$.

Here is a start up result.
\begin{theorem}
                                        \label{theorem 11.4.1}
For any $x_{0}\in\bR^{d}$ there exists a
$d$-dimensional Wiener process $w_{t}$
such that $(w_{t},\bar \cN^{x_{0}}_{t})$ is a Wiener process
and $P_{x_{0}}$-(a.s.) for all $t\geq0$
\begin{equation}
                                               \label{11.4.1}
x_{t}=x_{0}+\int_{0}^{t}\sqrt{a(x_{s})}\,dw_{s}
+\int_{0}^{t} b(x_{s}) \,ds.
\end{equation}
\end{theorem}

This theorem would be a simple consequence of Theorem 4.5.1
of \cite{SV_79} if $b$ were supposed to be bounded.
In our case Theorem \ref{theorem 11.4.1}
is a direct corollary of Lemma 1.6 of \cite{Kr_19_1}
and the following.

\begin{lemma}
                                        \label{lemma 11.4.1}
For any $x_{0}\in\bR^{d}$ and any
twice continuously differentiable function $u$
with compact support the process
\begin{equation}
                                               \label{11.6.2}
\xi_{t}:=u(x_{t})-\int_{0}^{t}Lu(x_{s})\,ds
\end{equation}
is a  martingale with respect to $\cN_{t}$
and measure $P_{x_{0}}$.
\end{lemma}

Proof. What we are given
and the Markov property imply that for any $0=t_{0}<t_{1}
<...<t_{n}=s<t$ and Borel bounded
$g(y_{0},...,y_{n})$ on $\bR^{d(n+1)}$ 
$$
E_{x_{0}}g(x_{t_{0}},...,x_{t_{n}})(\xi_{t}-\xi_{s})
$$
$$
=E_{x_{0}}g(x_{t_{0}},...,x_{t_{n}})\Big(E_{x_{s}}u(x_{t-s})
-u(x_{s})-E_{x_{s}}\int_{0}^{t-s}Lu(x_{r})\,dr\Big)=0.
$$
It follows that
$$
E_{x_{0}}(\xi_{t}-\xi_{s}\mid \bar\cN^{x_{0}}_{s})=0,
$$
which shows that $(\xi_{t}, \cN^{x_{0}}_{t })$,
$t\geq0$, is, indeed, a martingale.
 \qed

\mysection{Some results from \protect\cite{Kr_19}
and \protect\cite{Kr_19_1}}

                                              \label{section 12.1.1}
Theorem  \ref{theorem 11.4.1} allows us to use the results
from \cite{Kr_19}, \cite{Kr_19_1} some of which we list here.

We need  a part of
  Corollary 1.2   of \cite{Kr_19} which we state as follows. 

\begin{theorem}
                                    \label{theorem 11.6.1}
  For $0\leq s <t<\infty$, $x\in \bR^{d}$, and $n\geq0$ we have
\begin{equation}
                                                     \label{10.20.10}
 E_{x}\max_{r\in[s,t]}|x_{r}-x_{s}|^{2n}\leq N (t-s)^{n},
\end{equation}
where $N=N(n,d,\delta,\|b\|)$.

\end{theorem}

 Here are  Theorems 2.7 and 2.8 of \cite{Kr_19_1},
in which $d_{0}=d_{0}(d,\delta,\|b\|)\in(d/2,d)$ and for
Borel $\Gamma\subset\bR^{d}$
$$
\phi_{t}(\Gamma)=\int_{0}^{t}I_{\Gamma}(x_{s})\,ds.
$$

\begin{theorem}
                               \label{theorem 11.22.1}
Let $p\geq d_{0}$ and $R\in(0,\infty)$. Then
there exists $N$, depending only on
$p,d,\delta$, and $\|b\|$, such that for any 
 Borel nonnegative $f$ on $\bR^{d}$ we have
\begin{equation}
                                                     \label{9.5.4}
E_{0}\int_{0}^{\tau_{R}}f(x_{t})\,dt\leq
NR^{2-d/p}\|f\|_{L_{p}(B_{R})}.
\end{equation}

\end{theorem}
\begin{theorem}
                                        \label{theorem 11.23.1}
Let   $p\geq d_{0}$.
 Then there exist    constants $N$ and $\mu>0$,
depending only on $d,p$, and $\|b\|$,
and there exists $R_{0}=R_{0}(d,\|b\|)\geq2$, such that for any
$\lambda>0$, $R\in[0,\infty)$,
 and Borel nonnegative $f$ given on $\bR^{d}$
we have
\begin{equation}
                                                   \label{11.23.1}
E_{0}\int_{0}^{\infty}e^{-\lambda \phi_{t}(B_{R}^{c})}
f(x_{t}) \,dt   \le N (R\sqrt{\lambda}+R_{0}) ^{2-d/p} 
\lambda^{d/(2p)-1} 
\|\Psi_{R,\lambda}^{-1} f\| _{L_{p}(\bR^{d})},
\end{equation}
where $\Psi _{R,\lambda}(x)=\exp\big(\sqrt{\lambda}\,\mu\,
\dist(x,B_{R+R_{0}/\sqrt{\lambda}})\big)$.
\end{theorem}

Of the same spirit is the following particular case of
Theorem 4.7 of \cite{Kr_19}.

\begin{theorem}
                                \label{theorem 11.10.1}
 
 There exists   constants $N$ and $\mu>0$,
depending only on $d,\delta $,  and $\|b\|$, such that for any
$\lambda>0$ and Borel nonnegative $f$ given on $\bR^{d+1}$
we have
\begin{equation}
                                                   \label{7.29.02}
 E_{0}\int_{0}^{\infty}e^{-\lambda t}
f(t,x_{t}) \,dt   \le N
\lambda^{-d/(2d+2)}
\|\Psi_{\lambda}^{-1} f\| _{L_{d+1}(\bR^{d+1})},
\end{equation}
where $\Psi_{\lambda}(x)=\exp( \sqrt{\lambda}\mu|x|)$.

\end{theorem} 

Introduce
$$
R_{\lambda}f(x)=E_{x}\int_{0}^{\infty}e^{-\lambda t}
f(x_{t})\,dt.
$$

\begin{corollary}
                              \label{corollary 11.27.1}
Let $q\geq p\geq d_{0}$. Then there exists a  constant  $N$,
depending only on $d,p,q,\delta$, and $\|b\|$, such that for any
$\lambda>0$ 
 and Borel nonnegative $f$ given on $\bR^{d}$
we have
\begin{equation}
                                         \label{1.1.1}
   R_{\lambda} f   \leq N
\lambda^{ d/(2p)-1}
\|\Psi_{0,\lambda}^{-1} f\| _{L_{p}(\bR^{d})},
\end{equation}
\begin{equation}
                                         \label{11.27.5}
 \| R_{\lambda} f\| _{L_{q}(\bR^{d})}\leq N
\lambda^{-1+(d/2)(1/p-1/q)}
\|  f\| _{L_{p}(\bR^{d})}.
\end{equation}
\end{corollary}

Indeed, \eqref{1.1.1} is a particular case of \eqref{11.23.1}
when $R=0$,
which can be rewritten as
$$
(R_{\lambda}f(x))^{p}\leq N\lambda^{d/2-p}
\int_{\bR^{d}}e^{-\sqrt{\lambda}\mu|x-y|}f^{p}(y)\,dy.
$$
On the right we have the convolution of two functions.
Hence,
$$
\|R^{p}_{\lambda}f\|_{L_{q/p}(\bR^{d})}
\leq N\lambda^{d/2-p}
\Big(\int_{\bR^{d}}e^{-\sqrt{\lambda}\mu|x |q/p}\,dx
\Big)^{p/q}\| f^{p}\|_{L_{1}(\bR^{d})},
$$
which immediately yields \eqref{11.27.5}.

\begin{corollary}
                             \label{corollary 12.27.5}
Let $p\geq d_{0}$. Then for any $f\in L_{p}(\bR^{d})$
we have
\begin{equation}
                                         \label{12.27.7}
\lim_{\lambda\to\infty} 
\|\lambda R_{\lambda} f-f\| _{L_{p}(\bR^{d})}=0.
\end{equation}
\end{corollary}

Indeed, owing to Corollary \ref{corollary 11.27.1}
it suffices to prove \eqref{12.27.7} for
$f\in C^{\infty}_{0}(\bR^{d})$. For such $f$
by It\^o's formula (see Theorem 1.3 in \cite{Kr_19_1})
$$
E_{x}f(x_{t})e^{-t}=f(x)+E_{x}\int_{0}^{t}e^{-\lambda s}
(g-\lambda f)
(x_{s})\,ds,
$$
where $g=Lf \in L_{d_{0}}(\bR^{d})$. Hence
$$
\lambda R_{\lambda} f-f=R_{\lambda} g
$$
and owing to \eqref{11.27.5}
$$
\|\lambda R_{\lambda} f-f\| _{L_{p}(\bR^{d})}
\leq N\lambda^{-1+(d/2)(1/d_{0}-1/p)}\|g\|_{L_{d_{0}}
(\bR^{d})}.
$$
Here the exponent of $\lambda$ is strictly less than zero
because $d_{0}>d/2$ and this yields \eqref{12.27.7}.

We are also going to need
Corollary 4.3 of \cite{Kr_19_1}, which we state
as follows. 
\begin{theorem}
                                    \label{theorem 11.22.10}
For any $R\in(0,\infty)$ and
$\kappa\in(0,1)$ there exist constants $\mu\geq1$, $\theta>0$, and
$N$, depending only on $ d,\delta,\|b\|$, and $\kappa$,
 such that, for any 
$x\in   B_{\kappa R}$ and Borel set
$\Gamma\subset  B_{R}$
 \begin{equation}
                                          \label{11.5.4}
P_{x} (\phi_{\tau_{R}}(\Gamma) \geq\theta\gamma^{\mu} R^{2})
\geq N^{-1}\gamma^{2\mu},
\end{equation}
where $\gamma=  |\Gamma|/|B_{R}|$.
 \end{theorem}

Here is a specification
of Theorem 4.4 of \cite{Kr_19} to our case.
Recall that 
$$
C_{T,R}=[0,T)\times B_{R},\quad C_{R}=C_{R^{2},R}.
$$

\begin{theorem}
                                     \label{theorem 11.8.1}
For any $\kappa\in(0,1)$ there is a 
function $q(\gamma)$, $\gamma\in(0,1)$,
depending only on $d,\delta,\|b\|,\kappa$, and, naturally, on $\gamma$,
such that for any $R\in(0,\infty)$, $x\in B_{\kappa R}$,
and closed $\Gamma\subset C_{R^{2},R}$ satisfying
$|\Gamma|\geq \gamma|C_{R^{2},R}|$ we have
$$
P_{x}(\tau_{\Gamma}\leq \tau_{R^{2},R})\geq q(\gamma),
$$
where $\tau_{\Gamma}$ is the first time the process 
$(t,x_{t})$
hits $\Gamma$ and $\tau_{R^{2},R}$ is its first exit time from
$C_{R^{2},R}$. Furthermore, 
$q(\gamma)\to 1$ as $\gamma\uparrow 1$. 

\end{theorem}

Here is a corollary of estimate (2.19) of \cite{Kr_19}:
\begin{equation}
                                                     \label{12.5.7}
P_{0}(\tau_{R}\leq t)\leq 2\exp\Big(-\frac{\beta R^{2}}{t}\Big),
\end{equation}
where $\beta=\beta(d,\delta,\|b\|)>0$.

Here is a corollary of Theorem 4.17 of \cite{Kr_19_1}.
 \begin{theorem}
             \label{theorem 12.7.2}
Let $R\in(0,\infty)$, $ \kappa,\eta\in(0,1)$,
$x,y\in B_{\kappa R}$, 
 and $\eta^{-1} R^{2}\geq t\geq \eta R^{2}$.
Then there exist $N,\nu>0$, depending only on $ \kappa,\eta,d,\delta$,
and $\|b\|$, such that, for any
$\rho\in(0,1]$,
\begin{equation}
                                                 \label{12.9.1}
NP_{x}( x_{t}\in B_{\rho R}(y),\tau_{   R} > t )\geq \rho^{\nu }.
\end{equation}
 
\end{theorem}

\mysection{Strong Markov and strong Feller properties of $X$}

                                              \label{section 12.1.2}

Here is one of basic  results of this section.
Its proof would be greatly simplified
if we knew that $X$ is strong Markov. Then we would use
stopping times. Instead we use randomized stopping.

\begin{theorem}
                                    \label{theorem 11.22.2}
Let $p\geq d_{0}$, $\lambda>0$, and $f\in L_{p}(\bR^{d})$.
Then there exist $\alpha\in(0,1)$ and $N$,
depending only on $p,d,\delta$, and $\|b\|$, such that
\begin{equation}
                                                      \label{11.22.3}
|R_{\lambda}f(x)-R_{\lambda}f(y)|
\leq N \lambda^{(\alpha p+d)/(2p)-1}\|f\|_{L_{p}(\bR^{d})}|x-y|^{\alpha}.
\end{equation}
\end{theorem}

The following is an immediate and well-known
consequence of just continuity of $R_{\lambda}f(x)$  
with respect to $x$ (see, for instance,
Theorem I.8.11 in \cite{BG_68}).

\begin{corollary}
                                          \label{corollary 12.1.1.}
The process $(x_{t},\bar N_{t+},P_{x})$ is strong Markov.
\end{corollary}

Note that strong Markov processes 
$(x_{t},\bar N_{t},P_{x})$ corresponding to $L$
in case $b$ is bounded are constructed in \cite{Kr_73}
and strong Markov processes 
$(x_{t},\bar N_{t+},P_{x})$ corresponding to $L$
in case $b$ is bounded are constructed in
\cite{Ba_98}. We show that any {\em Markov\/} process
corresponding to $L$ with $b\in L_{d}$
is {\em strong\/} Markov with respect to $\bar N_{t+}$.
 Theorem \ref{theorem 11.22.2}
will be proved after some preparations.

The following technical result is obtained by
using Fubini's theorem for $f_{t}$ vanishing for 
$t\geq T$ and then letting $T\to\infty$.

\begin{lemma}
                                        \label{lemma 11.24.1}

(i) Let $\phi_{t},\Psi_{t},f_{t}$ be
nonnegative Borel functions on $[0,\infty)$
such that $\Psi_{0}=1$, $\Psi_{t}$ is continuous and decreasing.
Then
$$
\int_{0}^{\infty}e^{-\phi_{t}}f_{t}\,dt=-
\int_{0}^{\infty}e^{-\phi_{t}}\Big(\int_{t}^{\infty}
e^{-(\phi_{s}-\phi_{t})}f_{s}\,ds\Big)\,d\Psi_{t}
+\int_{0}^{\infty}e^{-\phi_{t}}\Psi_{t}f_{t}\,dt.
$$

(ii)  Let $ f_{t}$ be a
nonnegative Borel function  on $[0,\infty)$ and let 
$\phi_{t}$ and $\psi_{t}$ be absolutely continuous
on $[0,\infty)$ such that $\phi_{0}=\psi_{0}$. Assume that
$$
\int_{0}^{\infty}e^{-\phi_{t}}f_{t}
+\int_{0}^{\infty}e^{-\psi_{t}}|\phi_{t}'-\psi_{t}'|
\Big(\int_{t}^{\infty}e^{-(\phi_{s}-\phi_{t})}f_{s}\,ds\Big)\,dt <\infty.
$$
Then
$$
\int_{0}^{\infty}e^{-\phi_{t}}f_{t}\,dt=
\int_{0}^{\infty}e^{-\psi_{t}}f_{t}\,dt-
\int_{0}^{\infty}e^{-\psi_{t}}(\phi_{t}'-\psi_{t}')
\Big(\int_{t}^{\infty}e^{-(\phi_{s}-\phi_{t})}f_{s}\,ds\Big)\,dt.
$$
\end{lemma}

In the sequel we use the parameter $n$ that will be
ultimately send to infinity. Note
that by  $R _{n}(B_{R})$ one   usually means
the resolvent at $\lambda=n$ of the process killed
outside $B_{R}$. Our notation has a different meaning.
If $n\to\infty$, our $R _{n}(B_{R})$ converges
to the resolvent at $\lambda=0$ of the process killed
outside $B_{R}$.

\begin{lemma}
                                        \label{lemma 11.22.1}
Let $f $ be a nonnegative Borel bounded
function  on $\bR^{d}$. For Borel sets $\Gamma\subset\bR^{d}$
define
$$
\phi_{t}(\Gamma)=\int_{0}^{t}I_{\Gamma}(x_{s})\,ds,\quad
u(x)=E_{x}\int_{0}^{\infty}e^{-t}f(x_{t})\,dt.
$$
Also for $R > 0$ set
$$
R _{n}(B_{R})f(x)=
E_{x}\int_{0}^{\infty}e^{ - n\phi _{t}(B_{R}^{c})}
f(x_{t})\,dt.
$$
Then
\begin{equation}
                                                    \label{11.24.3}
u =R _{1}(B_{R})f -R _{1}(B_{R})I_{B_{R}}u .
\end{equation}
Furthermore, for any   $n\geq1$
\begin{equation}
                                             \label{11.24.31} 
 R _{1}(B_{R})f =R _{n}(B_{R})f +(n-1)R _{n}(B_{R})
I_{B_{R}^{c}}R _{1}(B_{R})f 
\end{equation}
Finally, for any Borel $\Gamma
\subset B_{R}$ and $h^{n}=nR_{n}(B_{R})f$
$$
h^{n}(x)=E_{x}\int_{0}^{\infty}e^{-n\phi _{t}(\Gamma\cup B_{R}^{c})}
nI_{\Gamma}(x_{t})h^{n}(x_{t})\,dt
$$
\begin{equation}
                                             \label{11.24.32}
+ nE_{x}\int_{0}^{\infty}e^{-n\phi _{t}(\Gamma\cup B_{R}^{c})}
 f(x_{t})\,dt.
\end{equation}

\end{lemma}

Proof. Equation \eqref{11.24.3} is obtained by using the Markov
property of $X$ and
by applying Lemma \ref{lemma 11.24.1} (ii) with
$\phi_{t}=t$, $\psi_{t}=\phi_{t}(B_{R}^{c})$. To check that
  Lemma \ref{lemma 11.24.1} (ii) is applicable we use
Theorem \ref{theorem 11.23.1} according to which  
$$
E_{x}\int_{0}^{\infty}e^{-\phi_{t}(B_{R}^{c})}\,dt<\infty.
$$

Equation \eqref{11.24.31} is obtained directly by using the Markov 
property and Lemma \ref{lemma 11.24.1} (i) with 
$\phi_{t}=\phi_{t}(B_{R}^{c})$,  $\Psi_{t}=
\exp(-(n-1)\phi_{t}(B_{R}^{c}))$.

Finally, to get \eqref{11.24.32} it suffices in Lemma \ref{lemma 11.24.1} (i)
to take 
$\phi_{t}=n\phi_{t}(B_{R}^{c})$,
and $\Psi_{t}=\exp(-n\phi_{t}(\Gamma))$. The lemma is proved.\qed
 
\begin{remark}
                                       \label{remark 1.12.3}
Send $n\to\infty$ in \eqref{11.24.31}
and \eqref{11.24.32} assuming that $f=0$
in $B_{R}$ in \eqref{11.24.32}. Then intuitively we
should get
$$
u(x):=E_{x}\int_{0}^{\tau_{R}}f(x_{t})\,dt=
R_{1}(B_{R})f-E_{x}R_{1}(B_{R})f(x_{\tau_{R}}),
$$
$$
h^{n}(x)\to h(x):=E_{x}f(x_{\tau_{R}}),
\quad h(x)=E_{x}\big(I_{\tau_{\Gamma}<\tau_{R}}
h(x_{\tau_{\Gamma}})+I_{\tau_{\Gamma}>\tau_{R}}
f(x_{\tau_{R}})\big),
$$
where $\tau_{\Gamma}$ is the first time $x_{t}$
hits $\Gamma$.
The formulas we get this way are true
indeed if we know that $X$ is strong Markov
(and that  $\tau_{\Gamma}$ makes sense and
is a stopping time). Of course, $u$ and $h$ are
the objects of main interest, but we cannot handle them
because we do not know yet that $X$ is strong Markov.
That is why instead of using stopping times
we use randomized ones when, for instance 
in case of \eqref{11.24.31}, we stop $x_{t}$
on each time interval $dt$ it spends outside $B_{R}$
with probability $(n-1)dt$ (provided 
it was not stopped before).
\end{remark}

The next two lemmas are aimed at partially justifying
what is said in Remark \ref{remark 1.12.3}.

\begin{lemma}
                                        \label{lemma 11.22.2} 
Let   $c$ be a nonnegative
function such that $c\geq1$ on $B^{c}_{\rho}$.
Then for any $n,\rho,\varepsilon>0$, and $|x|<\rho+\varepsilon$
$$
I_{n,\varepsilon}:= E_{x}\int_{\tau_{\rho+2\varepsilon}}^{\infty}
e^{ -n\phi_{t} }
nc(x_{t}) \,dt 
\leq 2e^{-\sqrt{n}\varepsilon/N},
$$
where $N=N(d,\delta,\|b\|)$ and
$$
\phi_{t}=\int_{0}^{t}c(x_{s})\,ds.
$$

\end{lemma}

Proof. We look at the representation of $x_{t}$
as a solution of \eqref{11.4.1} with $x_{0}=x$.
Then after defining $\gamma$ as the first time after
$\tau_{\rho+\varepsilon}$ the process $x_{t}$ exits
from $B_{\varepsilon}(x_{\tau_{\rho+\varepsilon}})$
 we note that
$$
I_{n,\varepsilon}\leq E_{x}e^{ -n\phi_{\tau_{\rho+2\varepsilon}}}
 \leq
E_{x}E\Big(e^{ -n(\phi_{\gamma}-\phi_{\tau_{\rho+\varepsilon}})}
\mid \bar N_{\tau_{\rho+\varepsilon}}\Big).
$$
Here, by the conditional version
of Theorem 2.10 of \cite{Kr_19},
 the conditional expectation (a.s.) is dominated by 
$2e^{-\sqrt{n}\,\varepsilon/N}$,
and this proves the lemma. \qed

\begin{lemma}
                                        \label{lemma 11.25.1}
 For Borel $\Gamma $ and $\rho>0$ define $\Gamma_{\rho}
=\Gamma\cap B_{\rho}$. Assume that $\Gamma$ and $\rho$ are such that
$|\Gamma_{\rho}|\geq (1/2)|B_{\rho}|$ and $\Gamma_{3\rho}
= \Gamma_{\rho}\cup \{2\rho\leq |x|< 3\rho\}$. Then there are   constant
$N,\nu>0$, depending only on $d,\delta$, and $\|b\|$,
such that for any  $n>0$ and   $|x|\leq (5/2)\rho$
\begin{equation}
                                                    \label{11.25.1}
1\geq nE_{x}\int_{0}^{\tau_{3\rho}} I_{\Gamma_{3\rho}}(x_{t})
e^{-n\phi_{t}(\Gamma_{3\rho})}\,dt\geq 1-2e^{-\sqrt{n}\rho/N},
\end{equation}
and for $|x|\leq \rho$ 
\begin{equation}
                                                    \label{11.25.11}
I:=E_{x}\int_{0}^{\tau_{ 3\rho}}nI_{\Gamma_{ \rho}}(x_{t})
e^{-n\phi_{t}(\Gamma_{3\rho})}\,dt\geq \nu-Ne^{-n\rho^{2}/N}.
\end{equation}
\end{lemma}

Proof. Observe that
$$
I\geq E_{x}\int_{0}^{\tau_{ 2\rho}}nI_{\Gamma_{ \rho}}(x_{t})
e^{-n\phi_{t}(\Gamma_{ \rho})}\,dt=
1-E_{x}e^{-n\phi_{\tau_{ 2\rho}}(\Gamma_{ \rho})},
$$
where for $\gamma=|\Gamma_{\rho}|/|B_{2\rho}|$ ($\geq 2^{-d-1}$)
owing to Theorem \ref{theorem 11.22.10}
$$
E_{x}e^{-n\phi_{\tau_{ 2\rho}}(\Gamma_{ \rho})}\leq
P\big( \phi_{\tau_{ 2\rho}}(\Gamma_{ \rho})\leq \theta
\gamma^{\mu}(2\rho)^{2}\big)
$$
$$
+e^{-n\theta\gamma^{\mu}(2\rho)^{2}}P\big( \phi_{\tau_{ 2\rho}}(\Gamma_{ \rho})
\geq \theta\gamma^{\mu}(2\rho)^{2}\big)
$$
$$
=1-(1-e^{-n\rho^{2}/N})P\big( \phi_{\tau_{ 2\rho}}(\Gamma_{ \rho})
\geq \theta\gamma^{\mu}(2\rho)^{2}\big)\leq 
1-(1-e^{-n\rho^{2}/N})N^{-1}\gamma^{2\mu}.
$$
This proves \eqref{11.25.11}.

To prove \eqref{11.25.1} denote by $\gamma$
the first exit time of $x_{t}$ after $\tau_{(5/2)\rho}$
from $B_{(1/2)\rho}(x_{\tau_{(5/2)\rho}})$ and observe that
$$
0\leq 1-nE_{x}\int_{0}^{\tau_{3\rho}} I_{\Gamma_{3\rho}}(x_{t})
e^{-n\phi_{t}(\Gamma_{3\rho})}\,dt=
E_{x}e^{-n\phi_{\tau_{3\rho}}(\Gamma_{3\rho})}
$$
$$
\leq E_{x}E\Big(e^{-n(\phi_{\gamma}(\Gamma_{3\rho})-\phi_{\tau_{(5/2)\rho}}(\Gamma_{3\rho}))}
\mid \cF_{\tau_{(5/2)\rho}}\Big).
$$
Here, by the conditional version
of Theorem 2.10 of \cite{Kr_19},
 the conditional expectation (a.s.) is dominated by
$2e^{-\sqrt{n}\,\rho/N}$,
and this proves the lemma. \qed

Before doing the next almost final step 
in our preparation to prove Theorem \ref{theorem 11.22.2},
take  $n\geq 1, R\in (0, 1]$, a bounded Borel
$g\geq0$, such that $g(x)=0$ for $|x|<R$ and introduce
$$
h^{n}=nR_{n}(B_{R})g.
$$
Observe that
$$
h^{n}\leq (\sup g) nE_{x}\int_{0}^{\infty}e^{-n\phi_{t}(B^{c}_{R})}
I_{B^{c}_{R}}(x_{t})\,dt\leq \sup g.
$$

The proof of the following lemma, actually,
is just a simple adaptation
of what is usually done in the theory
of elliptic equations when they
prove the H\"older continuity
of harmonic functions associated
with elliptic operators. Again
lacking the strong Markov property
and knowing nothing about the sets
$\Gamma^{n}$, introduced below,
apart from the fact that they are Borel,
forces us to use randomized stopping times.
Somewhat cleaner this adaptation
is seen in our Section \ref{section 12.1.3}.
Concerning the origin of our arguments
see Section 9.6 in \cite{Kr_18}.
 
\begin{lemma}
                                                 \label{lemma 11.24.3}
There exist constants $\alpha\in(0,1)$ and $N$, depending only on $d,
\delta$, and $\|b\|$, and there exists a constant
$N'$, depending only on $d,
\delta$, $\|b\|$, and  $\sup g$, such
 that for $|x|\leq R/12$ and any $n\geq1$ we have
\begin{equation}
                                                      \label{11.27.1}
|h^{n}(x)-h^{n}(0)|\leq N(|x|/R)^{\alpha}\sup_{B_{R}}h^{n}
+N'(e^{-n|x|^{2}/N}+e^{-\sqrt{n}\,|x|/N}).
\end{equation}
\end{lemma}

Proof. For $\rho>0$ introduce the notation
$$
\osc_{B_{\rho}}u=\sup_{B_{\rho}}u-\inf_{B_{\rho}}u,\quad
M^{n}_{\rho}=\sup_{B_{\rho}}h^{n},\quad 
m^{n}_{\rho}=\inf_{B_{\rho}}h^{n},\quad \mu^{n}_{\rho}=
(M^{n}_{\rho}+m^{n}_{\rho})/2.
$$
Take $\rho >0$, such that $ \rho \leq R/4$,
 and consider two cases

(a) $|B_{\rho}\cap\{h^{n}\geq \mu^{n}_{\rho}\}|
\geq (1/2)|B_{\rho}|$,

(b) $|B_{\rho}\cap\{h^{n}\leq \mu^{n}_{\rho}\}|
\geq (1/2)|B_{\rho}|$.

In case (a) introduce   $\Gamma^{n}=\big(B_{\rho}\cap
\{h^{n}\geq \mu^{n}_{\rho} \}\big)
\cup (B^{c}_{ 2\rho}\cap B_{R})$.
By using \eqref{11.24.32} and Lemma \ref{lemma 11.22.2}
for $|x|\leq\rho$ we find that
\begin{equation}
                                                   \label{11.26.1}
h^{n}(x)=E_{x}\int_{0}^{\tau_{3\rho}}e^{-n\phi _{t}(\Gamma^{n} )}
nI_{\Gamma^{n}}(x_{t})h^{n}(x_{t})\,dt+\xi^{n}(x)
=:h^{n}_{0}(x)+\xi^{n}(x),
\end{equation}
where $|\xi^{n}(x)|\leq N'e^{-\sqrt{n}\rho/N}$.
Furthermore, $h^{n}\geq \mu^{n}_{\rho}$ on $\Gamma^{n}_{\rho}$
and $h^{n}\geq m^{n}_{3\rho}$ in $B_{3\rho}$. Hence
$$
h^{n}_{0}(x)\geq (\mu^{n}_{\rho}-m^{n}_{3\rho})
E_{x}\int_{0}^{\tau_{3\rho}}e^{-n\phi _{t}(\Gamma^{n} )}
nI_{\Gamma^{n}_{\rho}}(x_{t}) \,dt
$$
$$
+m^{n}_{3\rho}E_{x}\int_{0}^{\tau_{3\rho}}e^{-n\phi _{t}(\Gamma^{n} )}
nI_{\Gamma^{n}}(x_{t}) \,dt.
$$
It follows by Lemma \ref{lemma 11.25.1}  
(also note that $\mu^{n}_{\rho}-m^{n}_{3\rho}\geq 0$) that
for $|x|\leq \rho$
$$
h^{n}_{0}(x)\geq \nu \mu^{n}_{\rho}+(1-\nu)m^{n}_{3\rho}
-N'(e^{-n\rho^{2}/N}+e^{-\sqrt{n}\rho/N}),
$$
which in light of the arbitrariness of $x$ and \eqref{11.26.1}
implies that
$$
(1-\nu/2)m_{\rho}^{n}\geq (\nu/2)M_{\rho}^{n}
+(1-\nu)m^{n}_{3\rho}-N'(e^{-n\rho^{2}/N}+e^{-\sqrt{n}\rho/N}).
$$
On the other hand, obviously
$$
(1-\nu/2)M_{\rho}^{n}\leq (\nu/2)M_{\rho}^{n}+(1-\nu)M^{n}_{3\rho}.
$$
By subtracting the last two inequalities we get
$$
(1-\nu/2)\osc_{B_{\rho}}h^{n}\leq(1-\nu)\osc_{B_{3\rho}}h^{n}
+N'(e^{-n\rho^{2}/N}+e^{-\sqrt{n}\rho/N}),
$$
\begin{equation}
                                                  \label{11.27.3}
\osc_{B_{\rho}}h^{n}\leq\theta\osc_{B_{3\rho}}h^{n}
+N'(e^{-n\rho^{2}/N}+e^{-\sqrt{n}\rho/N}),
\end{equation}
where $\theta=\theta(d,\delta,\|b\|)=(1-\nu)/
(1-\nu/2)<1$.

In case (b) introduce $\Gamma^{n}=\big(B_{\rho}\cap
\{h^{n}\leq \mu^{n}_{\rho} \}\big)
\cup (B^{c}_{ 2\rho}\cap B_{R})$. As in case (a),
we have \eqref{11.27.3}, where $h^{n}\leq \mu^{n}_{\rho}$
on $\Gamma^{n}_{\rho}$ and $h^{n}\leq M_{3\rho}^{n}$
om $B_{3\rho}$. Hence,
$$
h^{n}_{0}(x)\leq (\mu^{n}_{\rho}-M^{n}_{3\rho})
E_{x}\int_{0}^{\tau_{3\rho}}e^{-n\phi _{t}(\Gamma^{n} )}
nI_{\Gamma^{n}_{\rho}}(x_{t}) \,dt
$$
$$
+M^{n}_{3\rho}E_{x}\int_{0}^{\tau_{3\rho}}e^{-n\phi _{t}(\Gamma^{n} )}
nI_{\Gamma^{n}}(x_{t}) \,dt.
$$
Here $\mu^{n}_{\rho}-M^{n}_{3\rho}\leq0$ and the last expectation
is less than one. Then by \eqref{11.25.11}
$$
h^{n}_{0}(x)\leq (\mu^{n}_{\rho}-M^{n}_{3\rho})\nu+M^{n}_{3\rho}
+N'e^{-n\rho^{2}/N},
$$
which in light of the arbitrariness of $x$ and \eqref{11.26.1}
implies that
$$
(1-\nu/2)M_{\rho}^{n}\leq (\nu/2)m_{\rho}^{n}
+(1-\nu)M^{n}_{3\rho}+N'(e^{-n\rho^{2}/N}+e^{-\sqrt{n}\rho/N}). 
$$
On the other, hand obviously
$$
(1-\nu/2)m_{\rho}^{n}\geq (\nu/2)m_{\rho}^{n}+(1-\nu)m^{n}_{3\rho}.
$$
By subtracting the last two inequalities we get \eqref{11.27.3}
again.   

From     \eqref{11.27.3} we see that
$$
\theta^{-k}\osc_{B_{3^{-k} }}h^{n}\leq\theta^{-k+1}
\osc_{B_{3^{-k+1} }}h^{n}
+\theta^{-k}N'(e^{-n3^{-2k} /N}+e^{-\sqrt{n}\,3^{-k} /N}),
 $$
as long as    $3^{-k} \leq R/4$, that is 
$k\geq \lfloor\log_{3}(4/R)\rfloor=:k_{0}$. By
observing that, for instance, $\exp(-n3^{-2k}\rho^{2}/N)$
is an increasing function of $k$ we obtain, for $k\geq k_{0}$
$$
\theta^{-k}\osc_{B_{3^{-k}}}h^{n}\leq
\theta^{-k_{0}}\osc_{B_{3^{-k^{0}}}}h^{n}
+N'(e^{-n3^{-2k} /N}+e^{-\sqrt{n}\,3^{-k} /N})
\sum_{i=0}^{k-k_{0}-1}\theta^{-k+i}.
$$
For $|x|\leq R/12$ and $k=\lfloor\log_{3}(1/|x|)\rfloor$
 we have 
$$ |x|\leq 3^{-k}\leq R/4,\quad\theta^{-1}|x|^{-\alpha}\geq
\theta^{-k}\geq |x|^{-\alpha},
$$ 
where $\alpha=-\log_{3}\theta$. Furthermore,
$$
3^{-k}\geq |x|,\quad \theta^{-k_{0}}\leq 4^{\alpha}R^{-\alpha}.
$$
Now \eqref{11.27.1} follows. The lemma is proved. \qed

{\bf Proof of Theorem \ref{theorem 11.22.2}}. Self-similarity 
transformations show that we may
assume that $\lambda=1$. 
Furthermore, obviously we may assume that $f\geq0$.
Estimate \eqref{11.23.1} allows us to assume that $f$
is bounded and continuous. Then
take $R\in(0,1]$, $n\geq2$, and take $g=I_{B_{R}^{c}}R _{1}(B_{R})f$
in Lemma \ref{lemma 11.24.3}. Observe that by \eqref{11.24.31}
we have $h^{n}\leq 2R_{1}(B_{R})f$. Furthermore,   
 in light of \eqref{11.24.31},
  Lemma \ref{lemma 11.24.3}, and Theorem \ref{theorem 11.23.1},
 for $0<|x|\leq R/12$ we have 
$$
|R_{1}(B_{R})f(x)-R_{1}(B_{R})f(0)|\leq
N(R+R_{0}/\sqrt{n})^{2-d/p}\|f\|_{L_{p}(\bR^{d})}
$$
$$
+N(|x|/R)^{\alpha}(R+R_{0})^{2-d/p}\|f\|_{L_{p}(\bR^{d})}
+N'(e^{-n|x|^{2}/N}+e^{-\sqrt{n}\,|x|/N}),
$$
where $N'$ is independent of $n$. By sending $n\to\infty$
and taking onto account that $R\leq 1$ we come to
$$
|R_{1}(B_{R})f(x)-R_{1}(B_{R})f(0)|\leq
N R ^{2-d/p}\|f\|_{L_{p}(\bR^{d})}
$$
$$
+N(|x|/R)^{\alpha} \|f\|_{L_{p}(\bR^{d})}.
$$

By applying this result to \eqref{11.24.3}
and using Corollary \ref{corollary 11.27.1},
according to which $|u|\leq N\|f\|_{L_{p}(\bR^{d})}$,
and also using \eqref{1.1.1},   which implies
$R_{1}(B_{R})I_{B_{R}}\leq N R^{d/p}$ for $R\leq 1$,
 we obtain
$$
|R_{1}f(x)-R_{1}f(0)|\leq
N (R ^{2-d/p}+R^{d/p})\|f\|_{L_{p}(\bR^{d})}
$$
$$
+N(|x|/R)^{\alpha} \|f\|_{L_{p}(\bR^{d})}.
$$
If $2-d/p\leq d/p$, we take
 here $R=|x|^{\beta}$, where $\alpha\beta^{-1}=2-d/p+\alpha$, and
we get
 \begin{equation}
                                                   \label{11.27.6}
|R_{1}f(x)-R_{1}f(0)|\leq
N |x|^{\alpha(1-\beta)}\|f\|_{L_{p}(\bR^{d})},
\end{equation}
provided that $|x|\leq R/12=|x|^{\beta}/12$, that is
$|x|\leq \eta=\eta(d,\delta,p,\|b\|)$. For $|x|\geq \eta$
estimate \eqref{11.27.6} holds due to Theorem \ref{theorem 11.23.1}.
However, if $2-d/p>d/p$, we find $\beta$ from
$\beta(d/p+\alpha)=\alpha$ and again come to \eqref{11.27.6}. 

This proves \eqref{11.22.3} for $\lambda=1$ and $y=0$.
Shifting the coordinates take care of arbitrary $x,y$.
The theorem is proved. \qed

To prove that $X$ is strong Feller we need the following
generalization of a result of Lions \cite{Li_84},
proved in case $b$ is bounded and $p> d$, which was
generalized in \cite{FS_84} albeit when $b=0$ but with
$p\geq d_{0}$.
\begin{theorem}
                                          \label{theorem 12.2.1}
For any $p\geq d_{0}$ there are constants $N$ and $\mu>0$,
depending only on $d,p,\delta$, and $\|b\|$, such that
for any Borel nonnegative $f$ given 
on $\bR^{d}$ and $t>0$ we have
\begin{equation}
                                         \label{12.2.1}
 T_{t}f(0)\leq Nt^{-d/(2p)} 
\|\Phi_{t}f\|_{L_{p}(\bR^{d})},
\end{equation}
where $\Phi_{t}(x)=\exp(-\mu |x|/\sqrt{t})$.
\end{theorem}

Proof. The proof is quite similar to what is done
in \cite{Li_84}
and \cite{FS_84}. 
First fix $\varepsilon\in(0,1)$ and let
$F_{\varepsilon}$ be the set of Borel $f$
such that $\varepsilon\leq f(x)\leq \varepsilon^{-1}$
for any $x$. Then introduce
$$
K_{\varepsilon}=\sup\frac{T_{1}f(0)}{\|\Phi_{1}f\|_{L_{p}(\bR^{d})}},
$$
where the choice of $\mu$ (in $\Phi_{1}$)
will be specified later
and the supremum is taken over all $f\in F_{\varepsilon}$
 and over all 
diffusion processes
$X$, for which Assumption \ref{assumption 12.29.1}
 is satisfied
($\delta$ and $\|b\|$ are fixed). Obviously,
$K_{\varepsilon}<\infty$. Also observe that
self-similarity arguments easily show that
\eqref{12.2.1} holds with $N=K_{\varepsilon}$
if $f\in F_{\varepsilon}$. Shifting the origin
shows that for $f\in F_{\varepsilon}$, $t>0$,
and $x\in \bR^{d}$
\begin{equation}
                                                \label{12.3.1}
T_{t}f(x)\leq K_{\varepsilon}t^{-d/(2p)}
\Big(\int_{\bR^{d}}e^{-p\mu|x-y|/\sqrt{t}}f^{p}(y)
\,dy\Big)^{1/p}.
\end{equation}

Now, define $u(t,x)=T_{t}f(x)$ and observe
that by the Markov property for $s\in(0,1)$
$$
u(1,0)=T_{s}u(1-s,\cdot)(0)
=2\int_{0}^{1/2}T_{s}u(1-s,\cdot)(0)\,ds
$$
$$
=2E_{0}\int_{0}^{1/2}u(1-s,x_{s})\,ds
\leq NE_{0}\int_{0}^{\infty}I_{s\leq 1/2}e^{-s}
u(1-s,x_{s})\,ds.
$$

By using Theorem \ref{theorem 11.10.1} we obtain 
\begin{equation}
                                                \label{12.3.2}
u^{d+1}(1,0)\leq NIJ,
\end{equation}
where
$$
I=\sup_{[0,1/2]\times
\bR^{d}}\big( u^{d }(1-s,x )e^{-\mu|x|d}\big)
=\sup_{[1/2,1]\times
\bR^{d}}\big( u^{d }(t,x )e^{-\mu|x|d}\big),
$$
$$
J = \int_{\bR^{d}}e^{-\mu|x|}
\Big(\int_{1/2}^{1}u(t,x)\,dt\Big)\,dy .
$$
As is easy to see
$$
J \leq e\int_{\bR^{d}}e^{-\mu|x|}
R_{1}f(x)\,dx, 
$$
which by Theorem \ref{theorem 11.23.1} yields
$$
J \leq N\int_{\bR^{d}}e^{-\mu|x|}
\Big(\int_{\bR^{d}}e^{-\mu p|x-y|}f^{p}(y)\,dy\Big)^{1/p}\,dx, 
$$
where, perhaps, the second $\mu$ is different from the first
one. We allow ourselves to use $\mu$ as a generic constant
$>0$
depending only on $p, d,\delta$, and $\|b\|$. By H\"older's
inequality
$$
J\leq N\Big(\int_{\bR^{d}}e^{-\mu |x|}\,dx\Big)^{1/q}
\Big(\int_{\bR^{d}}g(y)f^{p}(y)\,dy\Big)^{1/p},
$$
where $q=p/(p-1)$ and
$$
g(y)=\int_{\bR^{d}}e^{-\mu (|x|+|x-y|)}\,dx.
$$
Since $|x|+|x-y|\geq(1/2)(|x|+|y|)$, we have
\begin{equation}
                                                \label{12.3.3}
J\leq N
\Big(\int_{\bR^{d}}e^{-p\mu|y|}f^{p}(y)\,dy\Big)^{1/p}.
\end{equation}

In what concerns $I$ observe that owing to \eqref{12.3.1}
for $t\in[1/2,1]$ we have
$$
u(t,x)e^{-\mu|x|}\leq NK_{\varepsilon}
\Big(\int_{\bR^{d}}e^{-\mu(|x|+|x-y|)}f^{p}(y)
\,dy\Big)^{1/p}
$$
$$
\leq NK_{\varepsilon}
\Big(\int_{\bR^{d}}e^{-p\mu|y|}f^{p}(y)\,dy\Big)^{1/p}. 
$$
Hence 
$$
I\leq NK_{\varepsilon}^{d}
\Big(\int_{\bR^{d}}e^{-p\mu|y|}f^{p}(y)\,dy\Big)^{d/p},
$$
and coming back to \eqref{12.3.2} we get
$$
u (1,0)\leq NK_{\varepsilon}^{d/(d+1)}
\Big(\int_{\bR^{d}}e^{-p\mu|y|}f^{p}(y)\,dy\Big)^{1/p}.
$$
Because of the definition of $K_{\varepsilon}$ it follows that
$$
K_{\varepsilon}\leq  NK_{\varepsilon}^{d/(d+1)},
\quad  K_{\varepsilon}\leq  N.
$$
After that it only remains to send $\varepsilon
\downarrow 0$ observing that the last $N$ as well as $\mu$
depend only on $p, d,\delta$, and $\|b\|$.
The theorem is proved. \qed

\begin{remark}
                                            \label{remark 12.3.1}
Once we know that \eqref{12.2.1} holds for Borel
nonnegative $f$, we can repeat the argument from the
beginning of the above proof and conclude that
for all Borel
nonnegative $f$, $t>0$,
and $x\in \bR^{d}$,
\begin{equation}
                                                \label{12.3.5}
T_{t}f(x)\leq Nt^{-d/(2p)}
\Big(\int_{\bR^{d}}e^{-p\mu|x-y|/\sqrt{t}}f^{p}(y)
\,dy\Big)^{1/p}.
\end{equation}
\end{remark}

\begin{corollary}
                                   \label{corollary 12.5.3}
For   $p\geq d_{0}$ such that $p>d/2+1$
there exists a constant $N=N(p,d,\delta,\|b\|)$
such that for any  $T\in(0,\infty)$ and nonnegative Borel $f(t,x)$
given on $[0,T]\times\bR^{d} $ we have
\begin{equation}
                                                \label{12.6.1}
I:=E_{0}\int_{0}^{T}f(t,x_{t})\,dt\leq NT^{(p-1)/p-d/(2p)}
\| \Phi _{T}f\|_{L_{p}([0,T]\times\bR^{d})}.
\end{equation}
\end{corollary}

Indeed,
\begin{equation}
                                                \label{12.5.5}
I=\int_{0}^{T}T_{t}f(t,\cdot)(0)\,dt
\leq N\int_{0}^{T}t^{-d/(2p)}
\Big(\int_{\bR^{d}} \Phi^{p}_{T}(y)f^{p}(t,y)
\,dy\Big)^{1/p}\,dt,
\end{equation}
and it only remains to use H\"older's inequality.

\begin{remark}
                                            \label{remark 12.5.5}
Observe that the usual parabolic Aleksandrov estimate
gives \eqref{12.6.1} with $p\geq d+1$. We were able
to reduce $p$ because $a$ is independent of $t$.
Also note that in \cite{Kr_74} there is an example showing
that \eqref{12.6.1} and \eqref{12.3.5} generally
(when $a$ is independent of $t$) fail  
to hold for any fixed $p<d$   
if $\delta$ can be chosen small enough. In that regard see
also the example in \cite{FS_84}, that appeared a few years
later than
\cite{Kr_74}.
\end{remark}

The following is deduced from \eqref{12.3.5} in the same way
as \eqref{11.27.5}
in Corollary \ref{corollary 11.27.1} is derived from
Theorem \ref{theorem 11.10.1}.

\begin{corollary}
                                \label{corollary 12.3.1}
For $q\geq p\geq d_{0}$ there exists a constant
 $N=N(p,q,d,\delta,\|b\|)$
such that for any $t>0$ and nonnegative Borel $f$
$$
\|T_{t}f\|_{L_{q}(\bR^{d})}\leq N
t^{(d/2)(1/q-1/p)}\| f\|_{L_{p}(\bR^{d})}.
$$
\end{corollary}

\begin{remark}
The fact that the semigroup $T_{t}$
is bounded in $L_{p}$, which follows from Corollary 
\ref{corollary 12.3.1} with $p=q$ , should not look
very surprising and follows by self-similarity
from the boundedness of $T_{1}$ (which, however,
is not trivial).
\end{remark}

\begin{corollary}
                                \label{corollary 12.28.1}
For $  p\geq d_{0}$ and $f\in L_{p}(\bR^{d})$
$$
\lim_{t\downarrow 0}\|T_{t}f-f\|_{L_{p}(\bR^{d})}=0.
$$
\end{corollary}

Indeed, in light of Corollary
\ref{corollary 12.3.1} (with $q=p$) we may concentrate
on $f\in C^{\infty}_{0}(\bR^{d})$. In that case
by It\^o's formula
$$
T_{t}f-f=\int_{0}^{t}T_{s}Lf\,ds,
$$
where $Lf\in L_{d_{0}}(\bR^{d})$. Hence, by 
Corollary
\ref{corollary 12.3.1}
$$
\|T_{t}f-f\|_{L_{p}(\bR^{d})}\leq
\int_{0}^{t}\|T_{s}Lf\|_{L_{p}(\bR^{d})}\,ds
\leq
\int_{0}^{t}s^{(d/2)(1/p-1/d_{0})}\,ds\| Lf\|_{L_{d_{0}}(\bR^{d})}
$$
and the last integral tends to zero as
$t\downarrow0$ since $d_{0}>d/2$.

\begin{theorem}
                                          \label{theorem 12.2.11}
The process $X$ is strong Feller in the sense that
for any bounded Borel $f$ given on $\bR^{d}$ the function
$T_{t}f(x)$ is continuous on $(0,\infty)\times\bR^{d}$.
\end{theorem}

Proof. In light of Theorem \ref{theorem 12.2.1}
we may concentrate on smooth compactly supported $f$'s.
In that case we are going to prove that $T_{t}f(x)$
is continuous in $[0,\infty)\times \bR^{d}$.
Observe that for $t\geq s\geq0$ we have
$$
|T_{t}f(x)-T_{s}f(x)|=|T_{s}(T_{t-s}f-f)(x)|
=|E_{x}E_{x_{s}}[f(x_{t-s})-f(x_{0})]|
$$
$$
\leq \sup|\nabla f|\sup_{y}E_{y}|x_{t-s}- x_{0} |\leq N
\sup|\nabla f|\sqrt{t-s},
$$
where the last inequality follows from Theorem \ref{theorem 11.6.1}.
We see that $T_{t}f(x)$ is continuous in $[0,\infty)$
uniformly with respect to $x$.

Next,  for our $f$ Theorem \ref{theorem 11.6.1} easily
implies that $\lambda R_{\lambda}f\to f$ as $\lambda\to
\infty$ uniformly with respect to $x$. After that
the continuity of $T_{t}f(x)$ with respect to $x$
follows from the fact that
$$
T_{t}\lambda R_{\lambda}f=\lambda R_{\lambda}T_{t}f,
$$
where the right-hand sides are continuous in $x$
due to Theorem \ref{theorem 11.22.2} and the left-hand sides
converge uniformly in $x$ to $T_{t}f(x)$.
The theorem is proved. \qed

We finish the section 
by a version of parabolic Aleksandrov
estimates.

\begin{theorem}
                       \label{theorem 12.30.1}  
Let $D$ be a bounded domain in $\bR^{d}$,
 $p\geq d_{0}$, $p>d/2+1$,  $T\in(0,\infty)$ and let    
$u\in W^{1,2}_{p,\loc}((0,T)\times D)\cap C([0,T]\times\bar D)$.
Then  there is a constant $N$, depending
only on $d,\delta,\|b\|$, and $T$, such that
in $[0,T]\times D$ we have
\begin{equation}
                               \label{12.30.1}
u\leq N\|(\partial_{t}u+Lu)_{-}\|_{L_{p}
((0,T)\times D)}+\sup_{(\partial ((0,T)\times  D))\setminus
(\{0\}\times D)}u .
\end{equation}
\end{theorem}

Proof. In the same way as in the proof of Theorem 3.1
of \cite{Kr_19} we convince ourselves that
we may assume that
$u\in W^{1,2}_{p }((0,T)\times D)$, $D$ is smooth, and  $b$ is bounded.
In that case it suffices to prove
\eqref{12.30.1} for $u\in C^{1,2}([0,T]\times \bar D)$.
For such $u$, by It\^o's
formula
$P_{0}$-(a.s.) for all $t\leq T$
$$
u(t\wedge \tau,x_{t\wedge \tau})=u(0,0)+\int_{0}^{t\wedge \tau}(\partial_{t}u+Lu)(s,x_{s})
\,ds +\int_{0}^{t\wedge \tau}D_{i}u(s,x_{s})\sigma^{ik}(x_{s})
\,dw^{k}_{s}
$$  
and the stochastic integral is a  
martingale, where $\tau$ is the first exit time of
$(t,x_{t})$ from  $[0,T)\times   D$. By setting $t=T$
and taking expectations we get
$$
u(0,0)=-E_{0}\int_{0}^{\tau }(\partial_{t}u+Lu)(t,x_{t})
\,dt 
+E_{0}u( \tau,x_{\tau })
$$
$$
\leq E_{0}\int_{0}^{\tau }
(\partial_{t}u+Lu)_{-}(t,x_{t})\,dt
+\sup_{(\partial ((0,T)\times  D))\setminus
(\{0\}\times D)}u .
$$
After that it only remains to use 
Corollary \ref{corollary 12.5.3}.
The theorem is proved. \qed

\mysection{Estimating time spent in space-time sets 
of small measure}
                                  \label{section 12.29.2}
Here we present 
extensions to the case that $b\in L_{d}$  of
 probabilistic versions of some
PDE results found in \cite{KS_80}, \cite{Sa_80},  
\cite{Kr_18}. Recall the notation introduced before
Theorem \ref{theorem 11.8.1} and also
introduce
$$
C^{o}_{T,R}=(0,T)\times B_{R},\quad C^{o}_{T,R}(t,x)=(t,x)+C^{o}_{T,R},
\quad C^{o}_{R}(t,x)=C^{o}_{R^{2},R}(t,x),
  $$
 $C^{o}_{R}=C^{o}_{R}(0,0)$.
Fix 
 $$
q, \eta ,  \kappa\in(0,1).
$$

 For cylinders $Q= C^{o}_{ \rho}(t ,x ) $
 define       
$$
Q'=(t ,x )-
C^{o}_{\eta^{-1} \rho^{2},\rho},\quad Q''=\big
(t -\eta^{-1}\rho^{2},x \big)
+C^{o}_{\eta^{-1}\rho^{2}\kappa^{2},\rho\kappa},
 $$
$$
Q'_{+}=Q\cup Q'\cup \big(\{t\}\times B_{\rho}(x)\big).
$$

Imagine that the $t$-axis is pointed up vertically. Then 
 $Q'$ is   adjacent to $Q$ from below, the 
two cylinders have a common base, and along the $t$-axis
$Q'$ is $\eta^{-1}$ times longer than $Q$.  
 The cylinder
$Q''$ is obtained by  contracting $Q'$  
to    the
center  of  its  lower base   with the contraction factor  
$\kappa^{-2}$ for the $t$-axis and $\kappa^{-1}$ for the spatial
axes. 

\begin{remark}
                                                   \label{remark 2,14,1}
If  $Q= C^{o}_{ \rho}(t ,x )$, then
  the shortest distance
between $Q$ and $Q''$ along the $t$ axis is     
\begin{equation}
                                                          \label{2,14,1}
\eta^{-1}\rho^{2}-\eta^{-1}\rho^{2}\kappa^{2}=\eta^{-1}\rho^{2}
(1-\kappa^{2}),
\end{equation}
 which is bigger than
$2\rho^{2}$ if     
\begin{equation}
                                                    \label{12.21.7}
\kappa^{2}\leq1-2\eta.
 \end{equation}
\end{remark}

 Let $\Gamma$ be a measurable subset of $C_{1}$
and
 introduce   $\mathcal{B}=\cB(\Gamma,q)$ 
 as the family  of   
{\em open\/}
cylinders $Q$ of  type 
$ C^{o}_{ \rho}(t_{0},x_{0})$ such that   
$$
Q\subset  C_{1 } \quad\text{and}\quad
|Q\cap\Gamma|\ge q|Q|.
$$ 

Finally,  define 
$$
\Gamma''=\bigcup_{Q\in\mathcal{B}}Q'',\quad
\Gamma''_{\varepsilon}=\bigcup_{Q\in\mathcal{B}:|Q|\geq
\varepsilon}Q''.
 $$

 Observe that for $Q\in\mathcal{B}$ the set 
   $Q''$ is open. Hence, 
$\Gamma''$ is open and measurable. 
\begin{lemma}
                                                     \label{lem:4.1.6} 
If $|\Gamma|\le q|C_{1 }|$, then   
$$|\Gamma''|\ge\Big(1-\frac{1-q}{3^{d+1}}\Big)^{-1}
(1+\eta)^{-1}\kappa^{d+2}|\Gamma|
 $$
and there exists $\theta=\theta(d,q,\eta,\kappa)>1$
such that for any sufficiently small $\varepsilon>0$
there exists a closed $\Gamma_{\varepsilon}
\subset \Gamma''_{\varepsilon}$ such that
\begin{equation}
                                                 \label{12.21.3}
|\Gamma_{\varepsilon}|\geq\theta |\Gamma|.
\end{equation}
\end{lemma}

The first assertion of the lemma originated
in \cite{KS_80}, \cite{Sa_80},  is presented, for instance
as Lemma 9.3.6 in
\cite{Kr_18}. The second one is proved
in the same way as the second assertion
of Lemma 4.7 of \cite{Kr_19_1}.

\begin{lemma}
                                         \label{lemma 12.20.2}
Let $\kappa\in(0,1)$. Then there is a constant 
$q_{0}=q_{0}(\kappa,d,\delta,\|b\|)
\in(0,1)$ such that for any $R\in(0,\infty)$,
 Borel set $\Gamma\subset C_{R}$
satisfying $|\Gamma|\geq q_{0}|C_{R}|$,
and $x\in B_{\kappa R}$ we have
\begin{equation}
                                                    \label{12.20.3}
E_{x}\int_{0}^{\tau_{R}\wedge R^{2}}I_{\Gamma}(t,x_{t})
\,dt\geq \mu_{0} R^{2} ,
\end{equation}
where   $\mu_{0}=\mu_{0}(d,\delta,\|b\|,
\kappa)\in(0,1)$.

\end{lemma}

Proof. As usual we let $R=1$. Then observe that    
by Lemma 2.13 of \cite{Kr_19} we have
$E_{x}(\tau_{1}\wedge1)\geq \nu=\nu(d,\delta,\|b\|,\kappa)>0$.
By using Corollary \ref{corollary 12.5.3} we get that
$$
E_{x}(\tau_{1}\wedge1)-
E_{x}\int_{0}^{\tau_{1}\wedge 1}I_{\Gamma}(t,x_{t})
\,dt=E_{x}\int_{0}^{\tau_{1}\wedge 1}I_{C_{1}\setminus
\Gamma}(t,x_{t})
\,dt
$$
$$
\leq N(|C_{1}|-|\Gamma|)^{1/d_{0}}
\leq N(1-q_{0})^{1/d_{0}}\leq N(1-q_{0})^{1/d_{0}}
E_{x}(\tau_{1}\wedge1),
$$
where the constants $N$ depend only on $\kappa,d,\delta$,  
and $\|b\|$. We see how to choose
$q_{0}$ to satisfy \eqref{12.20.3}
with a $\mu_{0}=\mu_{0}(d,\delta,\|b\|,
\kappa)\in(0,1)$.
The lemma is proved. \qed

In Lemma \ref{lemma 12.20.3} by $q_{0}$ we mean the one from Lemma \ref{lemma 12.20.2}.

\begin{lemma}
                                         \label{lemma 12.20.3}
Take $Q= C^{o}_{ \rho}(s,y) $,
use the notation $Q',Q'',Q'_{+}$ introduced above,
assume \eqref{12.21.7},
and suppose that Borel $\Gamma\subset Q$ is such that
$|\Gamma|\geq q_{0}|Q|$. Then there is a constant
$\nu_{0}>0$, depending only on $\eta,\kappa,d,\delta,\|b\|$,
such that for any $(t_{0},x_{0})\in Q''$
\begin{equation}
                                                    \label{12.20.4}
E_{x_{0}}\int_{0}^{\tau }I_{\Gamma}(t,x_{t})
\,dt\geq \nu_{0}E_{x_{0}}\tau,
\end{equation}
where $\tau$ is the first exit time of $(t_{0}+t,x_{t})$
from $Q'_{+}$.
\end{lemma}

Proof. Thanks to \eqref{12.21.7} and Remark \ref{remark 2,14,1}
we have $s-t_{0}\in(2\rho^{2}, \eta^{-1}\rho^{2})$.
Also $|y-x_{0}|<\kappa\rho$. It follows by
 Theorem \ref{theorem 12.7.2} that
$$
P_{x_{0}}\big(\sup_{r\in[0,s-t_{0}]}|x_{r}-y|< \rho,
|x_{s-t_{0}}-y|< \kappa\rho\big)\geq \nu ,
$$
where $\nu=\nu(\kappa,\eta,d,\delta,
\|b\|)>0$.

Next.
for $\gamma$ defined as the first exit time of $(t_{0}+t,x_{t})$
from $Q'$ in light of Lemma \ref{lemma 12.20.2} we have
$$
E_{x_{0}}\int_{0}^{\tau }I_{\Gamma}(t,x_{t})
\,dt=E_{x_{0}}I_{\gamma>s-t_{0}}
\int_{\gamma}^{\tau }I_{\Gamma}(t,x_{t})
\,dt
$$ 
$$
\geq E_{x_{0}}I_{\gamma>s-t_{0},|x_{s-t_{0}}-y|< \kappa\rho}
E_{x_{s-t_{0}}}\int_{0}^{\tau }I_{\Gamma}(t,x_{t})
\,dt
$$
$$
\geq \mu_{0}\rho^{2}
P_{x_{0}}\big(\sup_{r\in[0,s-t_{0}]}|x_{r}-y|< \rho,
|x_{s-t_{0}}-y|< \kappa\rho\big)\geq \mu_{0}\nu\rho^{2}.
$$
On the other hand, the height of $Q'_{+}$ is $(1+\eta^{-1})
\rho^{2}$, so that $(t_{0}+t,x_{t})$ cannot spend
in $Q'_{+}$ more time than $(1+\eta^{-1})
\rho^{2}$. This proves the lemma. \qed

\begin{theorem}
                                             \label{theorem 12.21.1}
For any $\kappa\in(0,1)$ there exists $\gamma\in(0,1)$
and $N$,
depending only on $\kappa,d,\delta,\|b\|$, such that
for any $R\in(0,\infty)$, $q\in(0,1)$,   Borel
$\Gamma\subset C_{R}(R^{2},0)$ satisfying
$|\Gamma|\geq q |C_{R}(R^{2},0)|$, and $x\in B_{\kappa R}$
we have
\begin{equation}   
                                                        \label{10.1.1}
 G_{R}(\Gamma,x):= E_{x}\int_{0}^{\tau_{R}\wedge(2R^{2})}
I_{\Gamma}(t,x_{t})\,dt\geq N^{-1}q^{1/\gamma}R^{2}.
 \end{equation}

\end{theorem}

Proof. Self-similar transformations allow us to assume that $R=1$
and write $G (\Gamma,x)$ instead of $G_{1}(\Gamma,x)$.
Then find and fix $q_{0},\eta,\kappa\in(0,1)$, depending only on
$d,\delta,\|b\|$, 
such that \eqref{12.21.7} holds, $\theta$ from Lemma \ref{lem:4.1.6}
is strictly bigger than 1,
and \eqref{12.20.3} holds whenever $|\Gamma|\geq q_{0}|C_{R}|$.
Clearly we can find such $\kappa\in(0,1)$ which is larger
than the one in the statement of the theorem.

It is convenient to introduce a function $\mu(q)$
as the infimum of the left-hand sides of \eqref{10.1.1}
(with $R=1$)
over all Borel $\Gamma\subset C_{1}(1,0)$ satisfying $|\Gamma|
\geq q|C_{1}(1,0)|$ and over all $x\in B_{\kappa}$.
Observe that a combination of Lemma \ref{lemma 12.20.2}
and Theorem \ref{theorem 12.7.2}, as in the proof of Lemma
\ref{lemma 12.20.3}, leads to the conclusion that
there exists $q_{0},\mu_{0}\in(0,1)$, depending only on
$\eta,\kappa,d,\delta,\|b\|$, such that
$$
\mu(q)\geq \mu_{0}
$$
for $q\in [q_{0},1]$.

We will be comparing $\mu(q')$ and $\mu(q'')$
for $ 0 <q'<q''<1$ such that
\begin{equation}
                                               \label{1.5.1}
(1+\theta) q'\geq 2q''.
\end{equation}

We take a Borel
$\Gamma\subset C_{1}(1,0)$ satisfying
$|\Gamma|\geq q' |C_{1}(1,0)|$ and
in the construction before Lemma \ref{lem:4.1.6} we replace $C_{1}$
by $C_{1 }(1,0)$, keep all other notation, and
from the chosen $\Gamma,\kappa,\eta$, and $q_{0}$ (not $q'$) we  
 build up  the
sets $\Gamma _{\varepsilon}$ and take
$\varepsilon$ so small that \eqref{12.21.3} holds.
  There are  two cases:  

(i) $\big|\Gamma_{\varepsilon}\setminus C_{1 }(1,0)\big|\leq (q''-q')|C_{1 }|$,

(ii) $\big|\Gamma_{\varepsilon}\setminus C_{1 }(1,0)\big|> (q''-q')|C_{1 }|$.

{\em Case  (i )\/}. Our  goal is to show that
\begin{equation}
                                                 \label{7,26,5}
G (\Gamma,x)
\geq\min\big(\mu_{0}, \nu_{0}\mu(q'')\big),\quad|x|\leq \kappa,
\end{equation}
where $\nu_{0}$
  depends only on  $\kappa,\eta,d,\delta,\|b\|$.

Observe that, if $|\Gamma|\geq q_{0}|C_{1 }|$, by definition
$G (\Gamma,x)\geq \mu(q_{0})\geq\mu_{0}$ for $|x|\leq R$.
 Hence, we may assume that
$$
|\Gamma|< q_{0}|C_{1 }|.
$$
In that case define
$$
\hat\Gamma_{\varepsilon}=\Gamma_{\varepsilon}\cap
C^{o}_{1 }(1,0).
 $$

Notice that by definition and  Lemma \ref{lem:4.1.6}   
$$
q'|C_{1 }|\leq|\Gamma|\leq \theta^{-1}|\Gamma_{\varepsilon}|.
 $$
Moreover, by assumption   
$$
|\Gamma_{\varepsilon}|=\big|\Gamma_{\varepsilon}
\setminus C_{1 }(1,0)\big|+|\hat  \Gamma_{\varepsilon}|
\leq (q''-q')|C_{1 }|+|\hat  \Gamma_{\varepsilon}|.
 $$
Due to \eqref{1.5.1}, it follows that
$$
|\hat \Gamma_{\varepsilon}|\geq q''|C_{1 }|,
$$
so that      
$$
G(\hat \Gamma_{\varepsilon},x)\geq\mu(q'') ,\quad|x|\leq \kappa.
 $$

We now estimate $G( \Gamma ,x)$ from
below by means of $G(\hat \Gamma_{\varepsilon},x)$
using Lemma \ref{lemma 12.20.3}. Since $\Gamma_{\varepsilon}
\subset \Gamma''_{\varepsilon}$,
the closed set $\Gamma_{\varepsilon}$ is 
covered by the family $\{Q'':Q\in\mathcal{B},|Q|\geq
\varepsilon\}$. Then there is
finitely many $Q(1),...,Q(n)\in \mathcal{B}$ such that
$|Q(i)|\geq\varepsilon$, $i=1,...,n$, and
$$
\Gamma_{\varepsilon}
\subset \bigcup_{i=1}^{n}Q''(i)
=:\Pi_{\varepsilon}.
$$

Then for $(t,x)\in \Pi_{\varepsilon}$ define  
$i (t,x) $
as the first $i\in\{1,...,n\}$ for which
$(t,x)\in Q''(i)$. 
Also set $Q'_{+}(0)=C_{2,1}$
and $i(t,x)=0$ if $(t,x)\in\partial C_{2,1}$.
Now define recursively $\gamma^{0}=0$,
$\tau^{1}$ as the first time after $\gamma^{0}$ when $(t,x_{t})$ exits
from $C _{2,1}\setminus \Gamma _{\varepsilon}$,
$\gamma^{1}$ as the first  time after $\tau^{1}$
when $(t,x_{t})$ exits from $Q' _{+}(i(\tau^{1},x_{\tau^{1}}))$,
and generally, for $k=2,3,...$ define
$\tau^{k}$ as the first time after $\gamma^{k-1}$ when $(t,x_{t})$ exits
from $C_{2,1}\setminus \Gamma _{\varepsilon}$,
$\gamma^{k}$ as the first  time after $\tau^{k}$
when $(t,x_{t})$ exits from $Q'_{+}(i(\tau^{k},x_{\tau^{k}}))$.
It is easy to check that so defined
$\tau^{k}$ and $\gamma^{k}$ are stopping times
and, since $|Q(i)|\geq\varepsilon$ and the trajectories 
of $(t,x_{t})$ are continuous,
$\tau^{k}\uparrow \tau_{1}\wedge 2$ as $k\to\infty$.
Furthermore,
(a.s.) all the $\tau^{k}$'s equal $\tau_{1}\wedge 2$
 for all large $k$.

For a domain $Q\subset\bR^{d+1}$
we denote by $\gamma(s,Q)$ the first exit time
of $(s+t,x_{t})$ from $Q$ and by the strong Markov
property obtain
$$
G(\Gamma,x)\geq \sum_{k=1}^{\infty}E_{x}\int_{\tau^{k}}
^{\gamma^{k}}I_{\Gamma}(t,x_{t})\,dt
$$
$$
=\sum_{k=1}^{\infty}E_{x}E_{x_{\tau^{k}}}\int_{0}
^{\gamma(s,Q'_{+}(i))}I_{\Gamma}(s+t,x_{t})\,dt\Big|_{
i=i(\tau^{k},x_{\tau^{k}}),s=\tau^{k}}.
$$
We estimate the interior expectation from below
by Lemma \ref{lemma 12.20.3} and get
$$
G(\Gamma,x)\geq \nu_{0}
\sum_{k=1}^{\infty}E_{x}E_{x_{\tau^{k}}}\int_{0}
^{\gamma(s,Q'_{+}(i))}I_{\Pi_{\varepsilon}}(s+t,x_{t})
\,dt\Big|_{s=\tau^{k},
i=i(\tau^{k},x_{\tau^{k}})}
$$
$$
\geq \nu_{0}
\sum_{k=1}^{\infty}E_{x}E_{x_{\tau^{k}}}\int_{0}
^{\gamma(s,Q'_{+}(i))}I_{\Gamma_{\varepsilon}}(s+t,x_{t})
\,dt\Big|_{s=\tau^{k},
i=i(\tau^{k},x_{\tau^{k}})}
$$
$$
= \nu_{0}
\sum_{k=1}^{\infty}E_{x} \int_{\tau^{k}}
^{\gamma^{k}}I_{\Gamma_{\varepsilon}}( t,x_{t})\,dt=
\nu_{0}G(\Gamma_{\varepsilon},x)\geq
\nu_{0}G(\hat \Gamma_{\varepsilon},x)\geq\nu_{0}\mu(q'').
$$
This proves \eqref{7,26,5}.

{\em Case (ii)\/}. Here the goal is to prove that  
\begin{equation}
                                                   \label{7,27,1}
G(\Gamma,x)\geq \mu_{0}\nu\eta^{n}(q''-q')^{n},
\quad |x|\leq \kappa,
 \end{equation}
  where $\nu>0$ and $n\geq1$ depend 
only on $d,\delta,\|b\|,\eta$, and $\kappa$.

First we claim that for some $(t,x)\in\Gamma_{\varepsilon}$
it holds that $t<q'-q''+1$. Indeed, otherwise
$\Gamma_{\varepsilon}\setminus C_{1}(1,0)
\subset C_{q''-q',1}(q'-q''+1,0)$ and 
$|\Gamma_{\varepsilon}\setminus C_{1}(1,0)|\leq
(q''-q')|C_{1}|$. It follows that there is a cylinder
$$
Q=C^{o}_{ \rho}(s,y)\in\cB
$$
 such that $Q'$
contains points in the half-space $t<q'-q''+1$. 
Since $q'<q''$, we have $q'-q''+1<1$, and since $Q'$ is adjacent
to $Q\subset C_{1 }(1,0)$, this implies that the height of $Q'$
is at least $q''-q'$, that is,  
\begin{equation}
                                            \label{12.22.7}
\rho^{2}\eta^{-1}\geq q''-q',\quad\rho^{2}\geq\eta(q''-q').
\end{equation}
On the other hand, $Q\subset C_{1}(1,0)$, $s>1$, and $\rho< 1$.

Moreover, 
 by construction,  $|\Gamma\cap Q|\geq q_{0}|Q|$ and
by Lemma \ref{lemma 12.20.2}
$$
E_{x}\int_{0}^{\tau}I_{\Gamma}(s+t,x_{t})\,dt
\geq \mu_{0}\rho^{2}\geq \mu_{0}\eta(q''-q')
$$
if $|x-y|\leq
 \kappa\rho $, where $\tau$ is the first exit time of $(s+t,x_{t})$
from $C_{ \rho}(s,y)$. Now by Theorem \ref{theorem 12.7.2}
for $x\in B_{\kappa }$
$$
E_{x}\int_{0}^{\tau_{2,1}}I_{\Gamma}(t,x_{t})\,dt
\geq E_{x}I_{\tau_{1}>s,|x_{s}-y|\leq\kappa\rho}
E_{x_{s}}\int_{0}^{\tau}I_{\Gamma}(s+t,x_{t})\,dt
$$
$$
\geq \mu_{0}\eta(q''-q')P_{x}\big
(\tau_{1}>s,|x_{s}-y|\leq\kappa\rho\big)\geq N^{-1}\rho^{\nu}
\mu_{0}\eta(q''-q').
$$
This proves \eqref{7,27,1}.

By combining the two cases (i) and (ii) we conclude that
$$
G(\Gamma,x)\geq \min\big(\mu_{0}, \nu_{0}\mu(q''),
\mu_{0}\nu\eta^{n}(q''-q')^{n}\big),
\quad |x|\leq \kappa,
$$
and the arbitrariness of $\Gamma$ allows us to conclude that
\begin{equation}
                                               \label{7,27,5}
\mu(q')\geq \min\big(\mu_{0}, \nu_{0}\mu(q''),
\mu_{0}\nu\eta^{n}(q''-q')^{n}\big),
\end{equation}
whenever \eqref{1.5.1} holds. Observe that
\eqref{7,27,5} is identical to (9.3.10) of \cite{Kr_18}
and by literally repeating what is in \cite{Kr_18},
just replacing $\xi$ there with our $\theta$,
we come to \eqref{10.1.1}. The theorem is proved. \qed

The following three results are derived from Theorem
\ref{theorem 12.21.1} in the same way as similar results are derived
from Theorem 4.1 of \cite{Kr_19_1}.

\begin{corollary}
                               \label{corollary 10.11.1}
For any 
$\kappa\in(0,1)$ there exists 
$N=N(d,\delta,\|b\|,\kappa)$  such that, for any $R\in(0,\infty)$,
$x\in   B_{\kappa R}$, and closed set
$\Gamma\subset  C_{R}(R^{2},0)$, the probability that the process 
$(t,x_{t})$ with $x_{0}=x$ 
reaches $\Gamma$ before exiting from $C_{2R^{2},R}$
is greater than or equal to $N^{-1} (|\Gamma|/|C_{R}|)^{\mu-1/d_{0}}$:
\begin{equation}
                               \label{10.2.10}
P_{x}(\tau_{\Gamma} <\tau_{2R^{2},R} )
\geq N^{-1} (|\Gamma|/|C_{R}|)^{\mu-1/d_{0}},
\end{equation}
where $\tau_{\Gamma} $ is the first time $(t,x_{t})$
hits $\Gamma$, $\tau_{2R^{2},R}$ is the first exit time of
$(t,x_{t})$ from $C_{2R^{2},R}$,
 $\mu=1/\gamma$, and $\gamma$ is taken from Theorem \ref{theorem 12.21.1}.
\end{corollary}

\begin{corollary}
                      \label{corollary 10.1.1}
For any Borel nonnegative $f$ 
vanishing outside $C_{R}(R^{2},0)$ and $x\in B_{\kappa R}$
$$
\int_{C_{R}(R^{2},0)}f^{1/(2\mu)}(t,y)\,dydt\leq NR^{d-1/\mu}
\Big(E_{x}\int_{0}^{\tau_{R^{2},R} }f(t,x_{t})\,dt\Big)^{1/(2\mu)},
$$
where $N=N(d,\delta,\|b\|,\kappa)$.

\end{corollary}

\begin{theorem}
                       \label{theorem 10.6.10}
Let    $p\in[d_{0},\infty)$, 
$u\in W^{1,2}_{p,\loc}(C_{2,1})\cap C(\bar C_{2,1})$,
  and
$c\in L_{d_{0}}(C_{2,1})$ $c\geq0$. Then
$$
\Big(\int_{C_{1,1} (1,0)}|D^{2}u|^{1/(2\mu)} 
\,dxdt\Big)^{2\mu} \leq N \sup_{\partial' C_{2,1}}|u|
$$
\begin{equation}
                              \label{10.6.1}
 +N \Big(\int_{C_{2,1}}|\partial_{t}u+
L u-cu|^{p}\,dxdt \Big)^{1/p},
\end{equation}
where $\partial_{t}=\partial/\partial t$, $\partial' C_{2,1}=
 \partial
C_{2,1}\setminus(\{0\}\times B_{1})$,
$\mu$ is taken from Corollary \ref{corollary 10.11.1} 
and $N$ depends only on
$d,\delta,\|b\|,p$, and $ \|c\|_{L_{d_{0}}(C_{2,1})}$.
\end{theorem}

It is worth emphasizing that in \eqref{10.6.1}
the restriction on $p$ is $p\geq d_{0}$ and $d_{0}<d$.
If $a$ depended on $t$, $p$ would be $>d$.

Theorem \ref{theorem 10.6.10} is similar to Theorem 9.4.1
of \cite{Kr_18} and in the same way as Theorem 9.4.9
of \cite{Kr_18} is derived from it (by using a simple trick)
one derives from Theorem \ref{theorem 10.6.10} the following.

\begin{theorem}
                       \label{theorem 12.23.1}
Let    $p\in[d_{0},\infty)$, 
$u\in W^{1,2}_{p,\loc}(C_{ 1})\cap C(\bar C_{ 1})$,
  and
$c\in L_{d_{0}}(C_{ 1})$, $c\geq 0$. Then
$$
\Big(\int_{C_{1 } }|D^{2}u|^{1/(2\mu)} 
\,dxdt\Big)^{2\mu} \leq N \Big(\int_{C_{ 1}}|\partial_{t}u+
L u-cu|^{p}\,dxdt \Big)^{1/p}
 +N \sup_{\partial' C_{ 1}}|u|,
$$
where   $\partial' C_{ 1}=\partial 
C_{ 1}\setminus(\{0\}\times B_{1})$,
$\mu$ is taken from Corollary \ref{corollary 10.11.1} 
and $N$ depends only on
$d,\delta,\|b\|,p$, and $ \|c\|_{L_{d_{0}}(C_{ 1})}$.
\end{theorem}

  In the next section we will need the following.
\begin{theorem}
                                      \label{thm:4.1.10} 
Let $\kappa,\eta,\zeta, q\in(0,1)$,  $R\in(0,\infty)$,
$T\in[\eta R^{2},\eta^{-1}R^{2}]$,
and closed $\Gamma\subset C_{T,R}$ be such that
$|\Gamma\cap C_{\zeta T,R}((1-\zeta)T,0)|\geq q
|C_{\zeta T,R}|$. Then there exists
$p_{0}=p_{0}(\kappa,\eta,\zeta,q,d,\delta,\|b\|)>0$ such that,
for $(t_{0},x_{0})\in C_{(1-\zeta)T, \kappa R}$,
\begin{equation}
                                      \label{12.27.3}
P_{x_{0}}(\tau_{\Gamma}<\tau_{T,R})\geq p_{0},
\end{equation}
where $\tau_{\Gamma}$ is the first time $(t_{0}+t,x_{t})$
hits $\Gamma$ and $\tau_{T,R}$ is its
first exit time from $C_{T,R}$.
\end{theorem}

Proof. As usual assume that $R=1$.
Then observe that one can choose 
$\rho>0$ depending only on $d,\eta,\zeta$, and $q$ and one can find
$(t^{0},x^{0})\in C_{T,1}$ with $t^{0}\ge \rho^{2}
+(1-\zeta)T$
such that $C_{\rho}(t^{0}+\rho^{2},x^{0})\subset C_{T,1}$
and
 $|\Gamma\cap C_{\rho}(t^{0}+\rho^{2},x^{0})|
\geq \bar q|C_{\rho} |$, where
$\bar q>0$ depends only on $d,\eta,\zeta$, and $q$. Then
by Corollary \ref{corollary 10.11.1}, for $x\in\kappa
B_{\rho}(x^{0})$ the $P_{x}$-probability that the process
$(t^{0}+t,x_{t})$ will hit $\Gamma$ before exiting
from $C_{2\rho^{2},\rho}(t^{0} ,x^{0})$ is estimated from 
below by a strictly positive constant depending only
on $\kappa,\bar q,d,\delta,\|b\|$.
After that it only remains to invoke Theorem
\ref{theorem 12.7.2} recalling that $t^{0}\geq \rho^{2}
+(1-\zeta)T$.
The theorem is proved. \qed

 \mysection{Harnack inequality, H\"older continuity
of $X$-caloric functions, and some other results}

                                           \label{section 12.1.3}

Safonov, \cite{Sa_10}, considered the case
of the coefficients  of $L$ so
regular  that $X$-harmonic functions are sufficiently
smooth and gave the estimate of the H\"older norm of $X$-harmonic
functions and the estimate of the Harnack constant
for them {\em independent of the imposed regularity
of $L$\/} in terms of only $d,\delta,\|b\|$.
We emphasize that $\|b\|$ is a bound of the $L_{d}$-norm
of $b$. Our case is not covered by \cite{Sa_10},
since the origin of our $X$ is unknown
and it is unknown if and in which sense
it can be approximated by processes with regular coefficients.
At the same time
  some arguments here are quite close to those in
\cite{Sa_10} as well as to those in \cite{KS_80},
\cite{Sa_80}, \cite{Kr_18}.

 \begin{definition}
                                \label{definition 12.27.1}
If $Q$ is a set in $\bR^{d+1}=\{(t,x)
:t\in\bR,x\in\bR^{d}\}$
and $u$ is a bounded Borel function on $Q$,
we call it caloric (relative to the process $X$) if
for any $(s,y)$ and $T,R\in(0,\infty)$
such that   $\bar C_{T,R}(s,y)\subset Q$
and any $(t_{0},x_{0})\in C:=C_{T,R}(s,y)$ we have
$$
u(t_{0},x_{0})=E_{x_{0}} u(t_{0}+\tau_{C},x_{\tau_{C }}),
$$
where $\tau_{C}$ is the first exit time of $(t_{0}+
t,x_{t})$
from $C$.

If $D$ is a set in $\bR^{d} $
and $u$ is a bounded Borel function on $D$,
we call it harmonic (relative to the process $X$) if
for any $y$ and $R\in(0,\infty)$
such that   $\bar B_{R}(y)\subset D $
and any $x\in B_{R}(y)$ we have
$$
u( x)=E_{x}u( x_{\tau_{R}(y) }),
$$
where $\tau_{R}(y)$ is the first exit time of $x_{t}$
from $B_{R}(y)$.

\end{definition}

\begin{remark}
                                  \label{remark 12.27.1}
If $u$ is harmonic in $D$
and $\bar B_{R}(y)\subset D $ and $x\in B_{R}(y)$, then 
by using the Markov property of $X$ 
we find
$$
u( x)=E_{x}E\big(u( x_{\tau_{R}(y) })\mid \cN_{T}\big)
=E_{x}I_{\tau_{R}(y)\leq T}u( x_{\tau_{R} })
+E_{x}I_{\tau_{R}(y)> T}E_{x_{T}}
 u( x_{\tau_{R} (y) })
$$
$$
=E_{x}I_{\tau_{R}(y)\leq T}u( x_{\tau_{R}(y) })
+E_{x}I_{\tau_{R}(y)> T}u(x_{T})=
 E_{x}u( x_{\tau_{R}(y)\wedge T  })
$$
 which implies that $u$
  is a
caloric function in $\bR\times D$.  
Also if $u$ is caloric in $Q$,
$\bar C_{T,R}(s,y)\subset Q$
and   $(t_{0},x_{0})\in C:=C_{T,R}(s,y)$, then by the strong Markov
property for any stopping time $\tau\leq
\tau_{C} $ we have
$$
u(t_{0},x_{0})=E_{x_{0}}E\big( 
 u(t_{0}+\tau_{C},x_{\tau_{C }})\mid \cN_{\tau}\big)
 =E_{x_{0}}  u(t_{0}+\tau ,x_{\tau }).
$$

\end{remark}

Here is the statement of the Harnack inequality.

\begin{theorem}
                                                           \label{thm:4.2.1} 
Let $\theta>1$, let $R\in(0,\infty]$, and let 
$u$ be a nonnegative 
caloric function in $\bar C_{\theta R^{2},R}$.
 Then there exists a constant $N$, which depends 
only on $\theta$,
$\delta$, $\|b\|$,   and $d$, such that   
\begin{equation}
                                          \label{12.26.1}
u(R^{2},0)\le Nu(0,x)
\end{equation}
whenever $|x|\le R/2$. 
\end{theorem}

Proof. As usual  without loss of generality
we  may concentrate on $R=2$. Then the case of general $\theta>1$
is reduced to that of $\theta\geq 2$ by appropriate change of 
the time variable 
    $t\to\tau(t)$.  One more observation that 
 the best constant $N$ is obviously  decreasing  in $\theta$,
 allows us to restrict our attention to the 
case of $\theta=2$.
 
In case $R=\theta=2$,
 to exclude a trivial situation, additionally assume that 
$$
u(4,0)>0.
$$
 
  For  
$\kappa=1/2,\eta=1/2$, we  take  $N$ and $\nu$ from
Theorem \ref{theorem 12.7.2}, call this $N$ $N_{1}$,
 and,  having in mind
    Theorem \ref{theorem 11.8.1},
find $\gamma\in(0,1)$    close to 1 and $\varepsilon>0$
close to zero,   for which     
\begin{equation}
                                       \label{eq:4.2.3}
1-\varepsilon\geq q(\gamma)2^{-1}+\big[1-q(\gamma)  \big]2^{\nu}.
\end{equation}

  Next, for $r\in[0,1)$, introduce  
$$\mu(r)=u(4,\,0)(1-r)^{-\nu},\quad n(r)=\sup\{ u,\bar C_{ r}(4,0)\}
\quad\big(\bar C_{0}(4,0):= \big\{(4,0)\big\}\big),
$$
  and define
  $r_{0}$ as the greatest number in $r\in [0,1)$ satisfying 
$$
n(r)\geq \mu(r).
$$ 
 Such a number does exist because
  $n(0)=\mu(0)$, $\mu(r)\to\infty$
as $r\uparrow1$, and $n(r)$ is bounded, increasing,
and right continuous.  
Choose $(t^{ \varepsilon},x^{ \varepsilon })
\in \bar C_{ r_{0}}(4,0)$
such that $n(r_{0})\leq(1+\varepsilon)
u(t^{ \varepsilon },x^{ \varepsilon })$ and consider  
the cylinder  
$$
Q:=\Big\{ (t,x)\,:\,0\leq t-t^{ \varepsilon }
<\frac{(1-r_{0})^{2}}{4},
\quad|x-x^{  \varepsilon}|<\frac{1-r_{0}}{2}\Big\} .
$$

  As is easy to see   $\bar Q\subset \bar C_{ r_{1}}(4,0)$,
where $r_{1}=(1+r_{0})/2 $.  By the definition of $r_{0}$,
this implies   that    
$$
\sup_{\bar Q}u<\mu(r_{1})=u(4,0)\Big(\frac{1-r_{0}}{2}\Big)^{-\nu}
\leq 2^{\nu}n(r_{0}).
$$

We claim that  owing to this and 
 (\ref{eq:4.2.3}), 
\begin{equation}
                                                       \label{eq:4.2.4}
\big|Q\cap\big\{ u>n(r_{0})/2\big\}\big|\ge  (1- \gamma ) |Q|.
\end{equation}

 To argue by contradiction, assume \eqref{eq:4.2.4}
is false. Then
$$
\big|Q\cap\big\{ u\leq n(r_{0})/2\big\}\big|
>    \gamma   |Q|
$$
and there is a closed set $\Gamma\subset
Q\cap\big\{ u\leq n(r_{0})/2\big\}$ such that
$|\Gamma|> \gamma  |Q|$. 
Introduce $\tau_{\Gamma}$ as the first time the process
$(t^{\varepsilon}+t,x_{t})$ hits $\Gamma$
and $\tau_{Q}$ as the first time it exits from $Q$.
It follows by 
Remark \ref{remark 12.27.1} and
Theorem \ref{theorem 11.8.1} that
(note that $n(r_{0})/2\leq\sup_{\bar Q}u$)
$$
u(t^{\varepsilon},x^{\varepsilon})=E_{x^{\varepsilon}}
I_{\tau_{\Gamma}<
\tau_{Q}}u(t^{\varepsilon}+\tau_{\Gamma},
x_{\tau_{\Gamma}})
+E_{x^{\varepsilon}}I_{\tau_{\Gamma}
\geq\tau_{Q}}u(t^{\varepsilon}+\tau_{Q},
x_{\tau_{Q}})
$$
$$
\leq P_{x^{\varepsilon}}
(\tau_{\Gamma}
<\tau_{Q})n(r_{0})/2+
(1-P_{x^{\varepsilon}}
(\tau_{\Gamma}
<\tau_{Q}))\sup_{\bar Q}u
$$
$$
\leq q(\gamma)n(r_{0})/2+
(1-q(\gamma))\sup_{\bar Q}u
$$
$$
\leq q(\gamma)n(r_{0})/2+
(1-q(\gamma))2^{\nu}n(r_{0}).
$$
We  now have
$$
n(r_{0})\leq(1+\varepsilon)n(r_{0})
\big[q(\gamma)2^{-1}+(1-q(\gamma))2^{\nu}\big]
\leq (1-\varepsilon^{2})  n(r_{0}),
$$
which is impossible. This proves 
 (\ref{eq:4.2.4}).

  Next we  apply Theorem \ref{thm:4.1.10}
and get that 
$$
u(t^{\varepsilon},x)\ge p_{0}n(r_{0})2^{-1}
$$
if $|x-x^{\varepsilon}|\le (1-r_{0})4^{-1}$, where  
$p_{0}=p_{0}( d,\delta,\|b\|, \gamma )>0$. 
 After that it only remains to apply  
Theorem \ref{theorem 12.7.2}
to conclude that for $|x|\le1$ we have    
$$
u(0,x)\ge\frac{1}{2}p_{0}n(r_{0})N_{1}^{-1}
\Big(\frac{1-r_{0}}{4}\Big)^{\nu}\geq 2^{-2\nu-1}p_{0}
N_{1}^{-1} u(4,0).   
$$
 The theorem is proved.                          \qed

Since harmonic function are also caloric
we have the following.
\begin{corollary}
                                  \label{corollary 12.28.3}
Let $R\in(0,\infty]$ and let $u$ be a
nonnegative harmonic function
in $B_{2R}$. Then for any $x,y\in B_{R}$ we have
$u(x)\leq Nu(y)$, where $N=N(d,\delta,\|b\|)$.
\end{corollary}

{\color{black}
\begin{corollary}
Let $R\in(0,\infty]$ and let $u\in W^{2}_{d_{0}}(B_{2R})$ be a
nonnegative   function satisfying $Lu=0$ (a.e.)
in $B_{2R}$. Then for any $x,y\in B_{R}$ we have
$u(x)\leq Nu(y)$, where $N=N(d,\delta,\|b\|)$.
\end{corollary}

Indeed, by Theorem 1.3 of \cite{Kr_19_1} It\^o's
formula is applicable and it shows that
$u$ is harmonic in $B_{2R}$.

}

 In the next part of the section
we deal with H\"older norm estimates for 
harmonic functions and potentials.
If $z_{1}=(t_{1},  x_{1})$ and $z_{2}=(t_{2}, x_{2})$,
we define    
\begin{equation}
                                                   \label{eq:4.2.5}
\rho(z_{1}, z_{2})=|x_{1}-x_{2}|+|t_{1}-t_{2}|^{1/2}
\end{equation}
 and call  $\rho(z_{1},z_{2})$   the parabolic
distance between $z_{1}$ and $z_{2}$.

\begin{lemma}
                                        \label{lem:4.2.2} 
Let $R\in(0, \infty]$ and let $u$ be a caloric
function in $\bar C_{2R}$. Then there
exist constants $N$ 
and   
$$
\alpha_{0}\in(0,1),
$$ 
depending only on $\delta,d,\|b\|$,
  such that, for any $\alpha\in(0,\alpha_{0}]$ and 
$z_{1},z_{2}\in C_{ R}$, we have   
\begin{equation}
                                       \label{eq:4.2.6}
\big|u(z_{1})-u(z_{2})\big|\le 
NR^{-\alpha}\rho^{\alpha}(z_{1},z_{2})
\sup\big(|u|,\bar C_{ 2R}\big).
\end{equation}

Furthermore, $\sup(|u|,\bar C_{ 2R})$ in 
\eqref{eq:4.2.6} can be replaced by $\osc(u,\bar 
C_{ 2R})$, where we use the 
notation
$$
 \osc(g,\Gamma) =\osc_{\Gamma}g=\sup_{\Gamma}g-
\inf_{\Gamma}g.
$$

\end{lemma}

Proof.  The case that $R=\infty$
is obtained by passing to the limit
and the case $R\in(0,\infty)$ reduces to $R=1$
by using self-similarity. In that case for $r\in(0,2]$, set
$$
w(r)=\osc(u,\bar C_{r}),\quad 
m(r)=\inf_{\bar C_{r}}u,\quad M(r)=\sup_{\bar C_{r}}u,
$$
$$
\mu(r)=(1/2)
\big(m(r)+M(r)\big).
$$

  Take  $r\le1/2$ and  suppose that
$$
\big|C_{ 2r}\cap\big\{ u\le
\mu(r)\big\} 
\big|\ge (1/2)|C_{ 2r}|.  
$$
Then there is a closed $\Gamma\subset
C_{ 2r}\cap\big\{ u\le
\mu(r)\big\}$ such that
\begin{equation}
                                   \label{eq:4.2.7}
\big|C_{3 r^{2},r}(r^{2},0)\cap\Gamma 
\big|\ge (1/4)|C_{3 r^{2},r}|
\end{equation}

By Theorem \ref{thm:4.1.10} for any $(t_{0},x_{0})
\in \bar C_{r}$ we have
$$
P_{x_{0}}(\tau_{\Gamma}<\tau_{2r})\geq p_{0},
$$
where $p_{0}>0$ depends only on $\delta,d,\|b\|$,
$\tau_{\Gamma}$ is the first time $(t_{0}+t,x_{t})$
hits $\Gamma$, $\tau_{2r}$ is its first exit
time from $C_{2r}$.
Then by Remark \ref{remark 12.27.1}
for $\tau=\tau_{\Gamma}\wedge\tau_{2r}$
$$
u(t_{0},x_{0})= E_{x_{0}}  u(t_{0}+\tau ,x_{\tau }).
$$
$$
=E_{x_{0}}  u(t_{0}+\tau_{\Gamma} ,x_{\tau _{\Gamma}})
I_{\tau_{\Gamma}<\tau_{2r}}
+E_{x_{0}}  u(t_{0}+\tau_{2r} ,x_{\tau_{2r} })
I_{\tau_{\Gamma}\geq\tau_{2r}}
$$
$$
\leq \mu(r)p_{0}+M(2r)(1-p_{0})
$$
(we used that $\mu(r)\leq M(2r)$).
It follows  that     
$$
M(r)\le p_{0}\frac{1}{2}\big(m(r)
+M(r)\big)+(1-p_{0})M(2r),
$$
$$
\big(1-\frac{p_{0}}{2})M(r)\leq \frac{p_{0}}{2} m(r)+(1-p_{0})M(2r).
$$

Adding  to this   the obvious inequality 
$$
\big(\frac{p_{0}}{2}-1)m(r)\leq -\frac{p_{0}}{2} m(r)
+(p_{0}-1)m(2r),
$$
we  get  
\begin{equation} 
                                             \label{eq:4.2.8}
\big(1-\frac{p_{0}}{2}\big)w(r)\le(1-p_{0})w(2r),\quad
w(r)\le\varepsilon w (2r),
\end{equation} 
where $\varepsilon<1$,  $\varepsilon=
\varepsilon(d,K,\delta)$.  We may, certainly, assume that 
 $\varepsilon>1/2$.

 We have proved (\ref{eq:4.2.8}) assuming 
that  (\ref{eq:4.2.7}) is true. However if (\ref{eq:4.2.7}) 
is false, then  
$-u$ satisfies an inequality similar to \eqref{eq:4.2.7}
and  this leads to
 (\ref{eq:4.2.8}) again.

\label{iteration}
 Therefore,  $w(r)\le\varepsilon w(2r)$ for all $r\le1/2$.  
Iterations then yield    
$$
w(r)\le\varepsilon^{2}w(4r)\quad\text{for}\quad r\le1/4,..., 
w(r)\le\varepsilon^{n}w(2^{n}r)  \quad\text{for}\quad r\le2 ^{-n}.
$$
 If $r\le1/2$ and  we take  $n:=\lfloor-\log_{2}r\rfloor$,
 then  $r\le2^{-n}$ and   
$$
w(r)\le\varepsilon^{n}w(2^{n}r)\le
\varepsilon^{-1}r^{\alpha}w(1)\le2
\varepsilon^{-1}r^{\alpha}
\sup\big(|u|,\bar C_{1 }\big),
$$
where $\alpha=-\log_{2}\varepsilon\in(0,1)$.  
This provides an estimate of 
the oscillation of $u$ in any   $C_{r}$
with $r\le1/2$. The same  estimate
obviously holds for the oscillation of $u$ in  any   
  $ C_{ r}(t,x)\subset C_{ 2}$ as long as  $r\le1/2$.

Now  take  $z_{1}=(t_{1},x_{1}),z_{2}=(t_{2},x_{2})\in C_{1}$ such that
$r:=\rho(z_{1},z_{2})\le1/2$ and define 
$$
t=t_{1}\wedge t_{2},\quad x= (x_{1}+x_{2})/2.
  $$
 Then  we have $z_{i}\in\bar C_{ R}(t,x)$, $i=1,2$,
and     
\begin{align*}
\big|u(z_{1})-u(z_{2})\big|  \le &\, 2
\varepsilon^{-1}r^{\alpha}\sup\big(|u|,
 \bar C_{1 }(t,x)\big)\\ 
   \le &\, 2\varepsilon^{-1}
\rho^{\alpha}(z_{1},z_{2})\sup\big(|u|,\bar C_{ 2}\big).
\end{align*}

 In  the case that $\rho(z_{1},z_{2})\geq 1/2$
 we have 
\begin{align*}
\big|u(z_{1})-u(z_{2})\big|   \le &\, 2
 \sup\big(|u|, \bar C_{2} \big)\\ 
  \le &\, 2^{1+\alpha}\rho^{\alpha}(z_{1},z_{2})
\sup\big(|u|,\bar C_{ 2}\big).
\end{align*}

Thus,   $N=2^{1+\alpha}+2 
\varepsilon^{-1} $ in (\ref{eq:4.2.6})  is always
 a good choice
with $R=1$
and $\alpha=\alpha(\delta,d)$ found above. One can take any smaller
$\alpha$ as well since $\rho(z_{1},z_{2})\leq N(d)R$.
The lemma is proved.                               \qed
  
 \begin{theorem}
                                     \label{theorem 10.8.1}

Let $R\in(0,\infty)$, $p\geq d_{0}$, $p>d/2+1$,
 let $g$ be a Borel bounded
function on $\bar C_{2R}$ and $f\in L_{p}(C_{2R})$.
For $(t_{0},x_{0})\in C_{2R}$ define
$$
u(t_{0},x_{0})=E_{x_{0}}
\int_{0}^{\gamma_{2R}}f(t_{0}+t,x_{t})\,dt+
E_{x_{0}}g(t_{0}+\gamma_{2R},x_{\gamma_{2R}}),
$$
where $\gamma_{2R}$ is the first exit time of 
$(t_{0}+t,x_{t})$ from $C_{2R}$.
  Then there exists a constant $N$, which depends
only on $p,d,\|b\|$, and $\delta$, such that   
\begin{equation}
                                                  \label{eq:4.2.9}
\big|u(z_{1})-u(z_{2})\big|\le N\big(R^{-\alpha}
\rho^{\alpha}(z_{1},z_{2})\sup_{\bar C_{ 2R}}|g|
+R^{2-(d+2)/p}\Vert f
\Vert _{L_{p}(C_{ 2R})}\big)
 \end{equation}
for $z_{1}$, $z_{2}\in  C_{ R}$ and 
$\alpha\in(0,\alpha_{0}]$.    
 
\end{theorem}

Proof.                               
Parabolic scalings  allow us to only concentrate
on the case
that $R=1$. After that it only remains  to
observe that $h(t_{0},x_{0}):
=E_{x_{0}}g(t_{0}+\gamma_{2R},x_{\gamma_{2R}})$
is a caloric function, to which
Lemma \ref{lem:4.2.2} is applicable, and
$u(t_{0},x_{0})-h(t_{0},x_{0})$ is estimated
by Corollary \ref{corollary 12.5.3}.
The theorem is proved.

Here is a  version
of Theorem \ref{theorem 10.8.1} which sometimes
is slightly more convenient.
\begin{theorem}
                                    \label{theorem 12,14.2}
Under the conditions of Theorem \ref{theorem 10.8.1}
there exists a constant $N$, which depends
only on $p,d,\|b\|$, and $\delta$, such that    
\begin{equation}
                                           \label{12,14.6}
\big|u(z_{1})-u(z_{2})\big|\le NR^{-\beta}
\rho^{\beta}(z_{1},z_{2})
\big(\sup_{\bar C_{ 2R}}|u|
+R^{2-(d+2)/p}\Vert f
\Vert _{L_{p}(C_{ 2R})}\big)
\end{equation}
for $z_{1}$, $z_{2}\in  C_{R}$, where
$$
\beta=\frac{\alpha(2p-d-2)}{\alpha p+2p-d-2}
$$
and $\alpha=\alpha_{0}(\delta,d)$ is the constant 
from Theorem \ref{theorem 10.8.1}  (or Lemma \ref{lem:4.2.2}).
\end{theorem}

Proof. Fix $z_{1}$, $z_{2}\in  C_{R}$.
Since  there is the sup norm of $u$
on the right, it suffices to prove \eqref{12,14.6}  
assuming that 
$$
\xi:=\Big(\frac{R}{\rho(z_{1},z_{2})}\Big)^{\beta/\alpha}\geq4.  
$$
Then set
$$
\bar R=\xi \rho(z_{1},z_{2}).
$$

If $z_{i}=(t_{i},x_{i})$, $i=1,2$, without losing generality
we may assume that $t_{1}\leq t_{2}$. Then  for
\begin{equation}
                                         \label{6,4,1}
|x_{1}|+\bar R\leq 2R\quad \text{and}\quad
 t_{1}+\bar R^{2}\leq 4R^{2}
\end{equation}
we have
\begin{equation}
                                         \label{12,14.7}
z_{1},z_{2}\in\bar C_{\bar R/4}(z_{1})\subset
\bar C_{\bar R }(z_{1})\subset \bar C_{2R}.
\end{equation}
 Since $z_{1}\in\bar C_{R}$,
we have $|x_{1}|\leq R$ and $t_{1}\leq R^{2}$
and, for any of the inequalities \eqref{6,4,1} to go wrong,
we have to have $\bar R\geq R$, that is,    
$$   
\Big(\frac{R}{\rho(z_{1},z_{2})}\Big)^{\beta/\alpha-1}\geq1,
$$
which is only possible if $\rho(z_{1},z_{2})\geq R$
when \eqref{12,14.6} holds trivially with $N=2$. Therefore,
in what follows we assume \eqref{12,14.7} and that 
$\bar R\leq R$.

 Then by Theorem \ref{theorem 10.8.1}
applied to $C_{\bar R}(z_{1})$ in place of $C_{R}$
we obtain   
$$
\big|u(z_{1})-u(z_{2})\big|\le N\big(\bar R^{-\alpha}
\rho^{\alpha}(z_{1},z_{2})\sup_{C_{ 2 R}}|u|
+\bar  R^{2-(d+2)/p}\Vert f
\Vert _{L_{p}(C_{ 2R})}\big),
$$
where the right-hand side is transformed 
to that of \eqref{12,14.6} by simple arithmetics.
 The theorem is proved.
\qed

\begin{corollary}
                                      \label{corollary 1.6.3}
Let $R\in(0,\infty)$, $p\geq d_{0}$, $p>d/2+1$,
and let $u\in W^{1,2}_{p}(C_{2R})$. Define
$f=\partial_{t}u+Lu$. Then there exists a constant $N$, which depends
only on $p,d,\|b\|$, and $\delta$, such that    
\eqref{12,14.6} holds for $z_{1}$, $z_{2}\in  C_{R}$
with the same $\beta$ as in \eqref{12,14.6}.
\end{corollary}
 
To prove this it suffices to follow the path
laid down in the proof of Theorem \ref{theorem 12.30.1}.

In the time-homogeneous situation we have a similar result.
\begin{theorem}
                                 \label{theorem 12.28.5}
Let $R\in(0,\infty)$, $p\geq d_{0}$, 
 let $g$ be a Borel bounded
function on $\bar B_{2R}$ and $f\in L_{p}(B_{2R})$.
For $ x_{0} \in B_{2R}$ define
$$
u( x_{0})=E_{x_{0}}
\int_{0}^{\tau_{2R}}f( x_{t})\,dt+
E_{x_{0}}g( x_{\tau_{2R}}),
$$
(recall that $\tau_{2R}$ is the first exit time of 
$ x_{t} $ from $B_{2R}$).
  Then there exists a constant $N$, which depends
only on $p,d,\|b\|$, and $\delta$, such that   
\begin{equation}
                                                         \label{1.6.7}
\big|u(x_{1})-u(x_{2})\big| \le NR^{-\alpha}
 |x_{1}-x_{2}|^{\alpha}
\big(\sup_{\bar B_{ 2R}}|u|
+R^{2-d/p}\Vert f
\Vert _{L_{p}(B_{ 2R})}\big)
\end{equation}
for $x_{1}$, $x_{2}\in  B_{ R}$ and 
$\alpha=\alpha(d,\delta,\|b\|)\in(0,1)$.    
 
\end{theorem}

This theorem is proved in the same way as
Theorem \ref{theorem 12,14.2} by using the fact
that $h(x_{0}):=E_{x_{0}}g( x_{\tau_{2R}})$
is a caloric function, to which Lemma 
\ref{lem:4.2.2}  is applicable,  and $u-h$
admits an estimate by Theorem \ref{theorem 11.22.1}.

Similarly to Corollary \ref{corollary 1.6.3} we have
the following corollary of Theorem \ref{theorem 12.28.5}.
\begin{corollary}
                                      \label{corollary 1.6.4}
Let $R\in(0,\infty)$, $p\geq d_{0}$,  
and let $u\in W^{ 2}_{p}(B_{2R})$. Define
$f= Lu$. Then there exists a constant $N$, which depends
only on $p,d,\|b\|$, and $\delta$, such that    
\eqref{1.6.7} holds for $x_{1}$, $x_{2}\in  B_{R}$
with the same $\alpha$ as in \eqref{1.6.7}.
\end{corollary}

We finish the paper by proving a result
showing that the function
$u$ from Theorem \ref{theorem 12.28.5} is an
$L_{d_{0}}$-viscosity solution
of the equation $Lu=-f$ in $B_{2R}$.
\begin{theorem}
                                 \label{theorem 12.28.2}
Let $u$ be as in Theorem \ref{theorem 12.28.5}.
Then for any $\phi\in W^{2}_{d_{0}}(B_{2R})$
and any point $x_{0}\in B_{2R}$ at which $u-\phi$
has local maximum we have
\begin{equation}
                                      \label{12.28.8}
\lim_{\varepsilon\downarrow0}
\esssup_{B_{\varepsilon}(x_{0})}(L\phi+f)\geq0.
\end{equation}
\end{theorem}

The proof of this theorem is based on the following.
\begin{lemma}
                                  \label{lemma 12.28.4}
There is a constant $N=N(d,\delta,\|b\|)$ such that
for any $B_{r}(x)$ satisfying $\bar B_{r}(x)
\subset B_{2R}$
and $\phi\in W^{2}_{d_{0}}(B_{r}(x))$
we have on  $B_{r}(x)$ that
\begin{equation}
                                      \label{12.28.7}
u\leq \phi+Nr^{2-d/d_{0}}\|(L\phi+f)_{+}\|_{L_{d_{0}}
(B_{r}(x))}+\max_{\partial B_{r}(x)}(u-\phi)_{+}.
\end{equation}
\end{lemma}

Proof. For $x_{0}\in B_{r}(x)$ by strong Markov property,
with $\tau $ defined as the first exit time of $x_{t}$
from $B_{r}(x)$,
and It\^o's
formula    we have
$$
u(x_{0})=E_{x_{0}}\Big(\int_{0}^{\tau }f(x_{t})\,dt
+u(x_{\tau})\Big),
$$
$$
\phi(x_{0})=E_{x_{0}}\Big(\int_{0}^{\tau }(-L\phi)(x_{t})\,dt
+\phi(x_{\tau})\Big).
$$
Hence,
$$
u(x_{0})-\phi(x_{0})\leq
E_{x_{0}}\Big(\int_{0}^{\tau }(L\phi+f)_{+}(x_{t})\,dt
+(u-\phi)_{+}(x_{\tau})\Big)
$$
and \eqref{12.28.7} follows from Theorem 
\ref{theorem 11.22.1}. The lemma is proved. \qed

{\bf Proof of Theorem \ref{theorem 12.28.2}}.
Let $x_{0}\in B_{2R}$ be a point at which $u-\phi$
has local maximum.
 Then for   $\varepsilon
>0$ and all small $r>0$ for    
$$
\phi_{\varepsilon,r}( x)=\phi ( x)-\phi( x_{0})
+u( x_{0})+\varepsilon(
|x-x_{0}|^{2} - r^{2})
$$
 we have that 
$$
\max_{\partial B_{r}( x_{0})}(u -\phi_{\varepsilon,r})_{+}
=0.
$$
Hence, by Lemma \ref{lemma 12.28.4}
$$
   \varepsilon r^{2}= 
(u -\phi_{\varepsilon,r})( x_{0})
\leq N_{1}r^{2-d/d_{0}}
\|(L\phi_{\varepsilon}+f)_{+}\|_{L_{d_{0}}
(B_{r}(x_{0}))},
$$
where $\phi_{\varepsilon }=\phi+\varepsilon
  |x -x_{0}|^{2} $.  Here $(L\phi_{\varepsilon}+f)_{+} 
\leq (L\phi +f+2\varepsilon \tr a)_{+} +N\varepsilon r|b| $
and in light of H\"older's inequality
$$
\|(L\phi_{\varepsilon}+f)_{+}\|_{L_{d_{0}}
(B_{r}(x_{0}))}\leq
\|(L\phi +f+2\varepsilon \tr a)_{+}\|_{L_{d_{0}}
(B_{r}(x_{0}))}
$$
$$
+N_{2}\varepsilon r^{d/d_{0}}\|b\|_{L_{d}(B_{r}(x_{0}))}.
$$
Here the last term multiplied by $N_{1}r^{2-d/d_{0}}$
is smaller than $(1/2)\varepsilon r^{2}$ for 
all sufficiently small $r$ (depending on how fast
$\|b\|_{L_{d}(B_{r}(x_{0}))}\to0$).
Therefore, for such $r$
$$
  (1/2) \varepsilon r^{2}=  
(u -\phi_{\varepsilon,r})( x_{0})
\leq Nr^{2-d/d_{0}}
\|(L\phi +f+2\varepsilon \tr a)_{+}\|_{L_{d_{0}}
(B_{r}(x_{0}))} 
$$
$$
\leq Nr^{2}\esssup_{B_{r}(x_{0}))}  
(L\phi +f+2\varepsilon \tr a)_{+},
$$
$$
\esssup_{B_{r}(x_{0}))}  
(L\phi +f+2\varepsilon \tr a)>0,
$$
and the last relation implies \eqref{12.28.8}
after setting $r,\varepsilon\downarrow0$.
The theorem is proved. \qed

\begin{remark}
                                  \label{remark 12.29.1}
Let $D$ be a bounded domain in $\bR^{d}$,
$g$ be a Borel bounded function on $\bR^{d}$,
and $f\in L_{d_{0}}(D)$. Introduce
$$
u(x)=E_{x}\Big(\int_{0}^{\tau}f(x_{t})\,dt
+g(x_{\tau})\Big),
$$
where $\tau$ is the first exit time of $x_{t}$
from $D$. Then by the strong Markov property
for any domain $G\subset D$  
$$
u(x)=E_{x}\Big(\int_{0}^{\gamma}f(x_{t})\,dt
+u(x_{\gamma})\Big),
$$
where $\gamma$ is the first exit time of $x_{t}$
from $G$. Therefore, Theorem \ref{theorem 12.28.2}
implies that $u$ is a viscosity solution
of $Lu+f=0$ in $D$. It is H\"older continuous
in $D$ in light of Theorem \ref{theorem 12.28.5}.

Its boundary behavior can be investigated
by using, for instance,   Theorem 4.10
of \cite{Kr_19_1}, which says that if
$0\in\partial D$
and   for some constants $\rho,\gamma>0$ and 
any $r\in (0,\rho)$ we have $|B_{r}\cap D^{c}|\geq \gamma
|B_{r}|$, then there exists 
$\beta=\beta(d,\delta,\|b\|,\gamma)>0$
such that, for any nonnegative $h\in L_{d_{0}}(D)$ 
and $x\in D$,
\begin{equation}
                                                 \label{10.20.4}
 E_{x}\int_{0}^{\tau }h( x_{t})\,dt
\leq N|x|^{\beta}\|h\|_{L _{d_{0}}(D)},
\end{equation}
where   $N$ depends only on $d,\delta,\|b\|,
\gamma,\rho$, and the diameter
of $D$.

The reader can find numerous properties of $L_{p}$-viscosity
 solutions in elliptic and parabolic
settings in articles initiated by \cite{CCKS_96},
references to many of them can be found in \cite{Kr_18}.
\end{remark}

{\bf Acknowledgment}. The author thanks T. Yastrzhembskiy
for pointing out several mistakes and misprints in the
first draft of the paper.


\begin{thebibliography}{mm}

      \bibitem{AP_77} 
	S.V. Anulova and
	G. Pragarauskas,
	{\em Weak Markov solutions of stochastic equations\/},
	  Litovsk. Mat. Sb., Vol. 17 (1977), No. 2, 5--26,   in Russian;
	English translation: Lithuanian Math. J., Vol. 17 (1977), 
No. 2, 141--155.

\bibitem{Ba_98} R. Bass,   ``Diffusions and elliptic operators'',
            Probability and its Applications,
Springer-Verlag, New York,
1998.

\bibitem{BG_68} R.M. Blumenthal and R.K. Getoor,
``Markov processes and potential theory'',
Pure and Applied Mathematics, A Series of Monographs
and Textbooks,
Vol. 29, Academic Press, New York and London, 1968.

\bibitem{CCKS_96} L. Caffarelli, M. G. Crandall,
M. Kocan, and A. \'Swi{\c e}ch,,
{\em On viscosity
solutions of fully nonlinear equations with
 measurable ingredients\/},
 Comm.
Pure Appl. Math., Vol. 49 (1996), No. 4, 365--397.

\bibitem{Dy_63}
 E. B. Dynkin, ``Markov processes'', 
 Fizmatgiz, Moscow, 1963 in Russian:
English translation in
 Grundlehren Math. Wiss., Vols. 121, 122,
Springer-Verlag,  Berlin, 1965.

\bibitem{FS_84} E.B. Fabes and D.W. Stroock,
{\em The $L^{p}$-integrability of Green's functions
and fundamental solutions for elliptic
and parabolic equations\/}, Duke Math. J., Vol. 51
(1984), No. 4, 997--1016.

 \bibitem{Kr_73}  N.V. Krylov,
  {\em  On the selection of a Markov process from a system
of processes and the construction of quasi-diffusion
processes\/},  Izvestiya Akademii Nauk SSSR, seriya
matematicheskaya, Vol. 37 (1973), No. 3, 691--708 in
Russian; English translation in Math. USSR Izvestija,
  Vol. 7 (1973),  No. 3,   691--709.

\bibitem{Kr_74}  N.V. Krylov, {\em  Some estimates for the density of the distribution
  of a
 stochastic
 integral\/}, Izvestiya Akademii Nauk SSSR, seriya matematicheskaya,
Vol. 38 (1974), No. 1,  228--248 in Russian; English translation
in Math. USSR
Izvestija,   Vol. 8  (1974),  No. 1,   233--254.

\bibitem{Kr_77}  N.V. Krylov,  ``Controlled diffusion processes'', 
Nauka, Moscow,  1977 in Russian; English transl.~
  Springer,
1980.

\bibitem{Kr_95} N.V. Krylov, ``Introduction to the theory 
of diffusion processes'',   Amer.
Math. Soc., Providence, RI, 1995.

\bibitem{Kr_18} N.V. Krylov,
``Sobolev and viscosity solutions for fully nonlinear  elliptic 
and parabolic equations'', Mathematical Surveys and Monographs,
233, Amer.
Math. Soc., Providence, RI, 2018.

\bibitem{Kr_19} N.V. Krylov, {\em 
On stochastic It\^o processes with drift in $L_{d}$\/},\\
http://arxiv.org/abs/2001.03660

\bibitem{Kr_19_1} N.V. Krylov, {\em  On stochastic equations with drift in
$L_{d}$\/},\\ http://arxiv.org/abs/2001.04008

\bibitem{KS_80}  N.V. Krylov and  M.V. Safonov,  {\em
   A certain property of solutions of parabolic equations
 with measurable coefficients}, 
Izvestiya Akademii Nauk SSSR, seriya matematicheskaya,
Vol.
44  (1980),  No. 1,  161--175  in Russian; English translation in
 Math. USSR
Izvestija, Vol. 16  (1981), No. 1,  151--164.

\bibitem{Li_84} P.-L. Lions, {\em
 Some recent results in the optimal
control of diffusion processes\/},
 Stochastic analysis (Katata/Kyoto, 1982),
pp 333--367, North-Holland Math. Library, 32,  North-Holland,
Amsterdam, 1984.

\bibitem{Po_82} N. I. Portenko, 
``Generalized diffusion processes'',  Nauka, Moscow, 
1982 in Russian; English translation: Amer. Math. Soc.
Providence, Rhode Island, 1990.

\bibitem{Sa_80}  M. V. Safonov, {\em
Harnack inequalities for elliptic equations and H{\"o}lder continuity
of their solutions\/}, 
Zap. Nauchn. Sem. Leningrad. Otdel. Mat. Inst. Steklov (LOMI),
Vol. 96 (1980), 272--287 in Russian;
English transl.~in Journal of Soviet Mathematics,
  Vol. 21 (March 1983), No. 5, 851--863.

\bibitem{Sa_10} M.V. Safonov, {\em Non-divergence elliptic
 equations of second
order with unbounded drift\/}, Nonlinear partial differential equations and
related topics, 211--232, Amer. Math. Soc. Transl. Ser. 2, 229, 
Adv. Math. Sci.,
64, Amer. Math. Soc., Providence, RI, 2010.

 \bibitem{SV_79} D.W. Stroock and S.R.S. Varadhan
``Multidimensional diffusion processes'',
 Grundlehren Math. Wiss.,  Vol. 233,
  Springer-Verlag, Berlin and New York, 1979.

\bibitem{Ya_19}
T. Yastrzhembskiy,
{\em A note on the strong Feller property of diffusion processes\/},
preprint.

\bibitem{ZZ_20} Xicheng Zhang and Guohuan Zhao,
{\em Stochastic Lagrangian path for Leray solutions of 3d
Navier-Stokes equations\/}, preprint.


 
\end{thebibliography}
\end{document}